\documentclass[a4paper]{amsart}
\usepackage{amssymb} 
\usepackage[T1]{fontenc}
\usepackage[latin1]{inputenc}
\usepackage{amsfonts}
\usepackage{amsxtra}
\usepackage{ae}
\usepackage[all]{xy}
\usepackage{enumerate}

\include{diagram}

\newcommand*{\cstaralg}{\textsf{-}C^*\textsf{-Alg}}
\newcommand*{\ket}{\rangle}
\newcommand*{\bra}{\langle}

\newcommand*{\Smooth}{\mathfrak{Smooth}}
\newcommand*{\Analysis}{\mathfrak{dis}}
\newcommand*{\A}{\mathcal{A}}
\newcommand*{\C}{\mathcal{C}}
\newcommand*{\D}{\mathcal{D}}
\newcommand*{\E}{\mathcal{E}}
\newcommand*{\ayd}{\mathsf{AYD}}
\newcommand*{\AYD}{\mathsf{A}}
\newcommand*{\cotimes}{\hat{\otimes}}
\newcommand*{\SHom}{\mathfrak{Hom}}

\newcommand*{\holim}{\text{ho}\text{-}\varinjlim}

\DeclareMathOperator{\Aut}{Aut}
\DeclareMathOperator{\Hom}{Hom}

\DeclareMathOperator{\id}{id}
\DeclareMathOperator{\ch}{ch}
\DeclareMathOperator{\Tot}{Tot}
\DeclareMathOperator{\an}{tan}
\DeclareMathOperator{\ind}{ind}

\newenvironment{bnum}
{\begin{list}{}
    {\setlength{\labelwidth}{15pt}
     \setlength{\leftmargin}{\labelwidth}
    }
}
{\end{list}}

\numberwithin{equation}{section}
\newtheorem{theorem}{Theorem}[section]
\newtheorem{prop}[theorem]{Proposition}
\newtheorem{lemma}[theorem]{Lemma}
\newtheorem{cor}[theorem]{Corollary}
\newtheorem{definition}[theorem]{Definition}

\begin{document}

\title[Equivariant local cyclic homology]{Equivariant local cyclic homology and the equivariant Chern-Connes character}
\author{Christian Voigt}
\address{Institut for Mathematical Sciences\\
         University of Copenhagen\\
         Universitetsparken 5 \\
         2100 Copenhagen\\
         Denmark
}
\email{cvoigt@math.ku.dk}


\maketitle

\begin{abstract}
We define and study equivariant analytic and local cyclic homology for smooth actions 
of totally disconnected groups on bornological algebras. Our approach contains equivariant 
entire cyclic cohomology in the sense of Klimek, Kondracki and Lesniewski as a special case 
and provides an equivariant extension of the local cyclic theory developped by Puschnigg. 
As a main result we construct a multiplicative Chern-Connes character for 
equivariant $ KK $-theory with values in equivariant local cyclic homology.
\end{abstract}

\section{Introduction}

Cyclic homology can be viewed as an analogue of de Rham cohomology in the framework of noncommutative geometry \cite{Connes1}, \cite{Connes2}. 
In this framework geometric questions are studied by means of associative algebras which need not be commutative. 
An important feature of cyclic homology is the fact that the theory can easily be defined on a large class of algebras, 
including Fr\'echet algebras as well as algebras without additional structure. In many cases explicit calculations are possible 
using standard tools from homological algebra. The connection to de Rham theory is 
provided by a fundamental result due to Connes \cite{Connes1} showing that the periodic cyclic homology of 
the Fr\'echet algebra $ C^\infty(M) $ of smooth functions on a compact manifold $ M $ is equal to the de Rham 
cohomology of $ M $. \\
However, in general the theory does not yield good results for Banach algebras or $ C^* $-algebras. Most notably, the 
periodic cyclic cohomology of the algebra $ C(M) $ of continuous functions on a compact manifold $ M $ is different from de Rham 
cohomology. An intuitive explanation of this phenomenon is that $ C(M) $ only encodes the information of $ M $ as a topological space, 
whereas it is the differentiable structure that is needed to define de Rham cohomology. \\
Puschnigg introduced a variant of cyclic homology which behaves nicely 
on the category of $ C^* $-algebras \cite{Puschnigg4}. The resulting theory, called local cyclic homology, allows for the 
construction of a general Chern-Connes character for bivariant $ K $-theory. 
Using the machinery of local cyclic homology, Puschnigg proved the Kadison-Kaplansky idempotent conjecture for hyperbolic 
groups \cite{Puschnigg3}. Unfortunately, the construction of the local theory is quite involved. Already the objects 
for which the theory is defined in \cite{Puschnigg4}, inductive systems of nice Fr\'echet algebras, are rather complicated. \\
There is an alternative approach to local cyclic homology due to Meyer \cite{Meyerlocalcyclic}. Based on the theory of 
bornological vector spaces, some features of local cyclic homology become more transparent in this approach. 
It is known that bornological vector spaces provide a very natural framework to study analytic and entire cyclic 
cohomology \cite{Meyerthesis}. Originally, entire cyclic cohomology was introduced by Connes \cite{Connesentire} in order to define
the Chern character of $ \theta $-summable Fredholm modules. The analytic theory for bornological algebras contains 
entire cyclic cohomology as a special case. Moreover, from a conceptual point of view it is closely related to the local theory. 
Roughly speaking, the passage from analytic to local cyclic homology consists in the passage to a certain derived category. \\
A central concept introduced by Meyer in his account to the local theory is the notion of 
an isoradial subalgebra \cite{Meyerlocalcyclic}, \cite{Meyerborntop}. Correspondingly, one of his main results is that local cyclic homology is invariant 
under the passage to isoradial subalgebras. In fact, an inspection of the proof of this theorem already reveals the essential ideas behind the definition of 
the local theory. A basic example of an isoradial subalgebra is the 
inclusion of $ C^\infty(M) $ into $ C(M) $ for a compact manifold as above. In particular, the natural 
homomorphism $ C^\infty(M) \rightarrow C(M) $ induces 
an invertible element in the bivariant local cyclic homology group $ HL_*(C^\infty(M), C(M)) $.  
Hence, in contrast to periodic cyclic cohomology, the local theory does not distinguish between $ C(M) $ 
and $ C^\infty(M) $. Let us also remark that invariance under isoradial subalgebras is responsible for the 
nice homological properties of the local theory. \\
In this paper we define and study analytic and local cyclic homology in the equivariant setting.  
This is based on the general framework for equivariant cyclic homology 
developped in \cite{Voigtepch} and relies on the work of Meyer in the nonequivariant case. In particular, a large part of the 
necessary analytical considerations is already contained in \cite{Meyerborntop}. In addition some of the material from 
\cite{Meyerlocalcyclic} will be reproduced for the convenience of the reader. On the other hand, as far as homological algebra is concerned, the 
framework of exact categories used by Meyer is not appropriate in the equivariant situation. 
This is due to the fact that equivariant cyclic homology is constructed using paracomplexes \cite{Voigtepch}. \\
We should point out that we restrict ourselves to actions of totally disconnected groups in this paper. 
In fact, one meets certain technical difficulties in the construction of the local 
theory if one moves beyond totally disconnected groups. For simplicity we have thus avoided to consider a more general setting.  
Moreover, our original motivation to study equivariant local cyclic homology and the equivariant Chern-Connes character 
comes from totally disconnected groups anyway. \\
Let us now explain how the paper is organized. In section \ref{sectotdis} we review some facts 
about smooth representation of totally disconnected groups and anti-Yetter-Drinfeld modules. 
These concepts are basic ingredients in the construction of equivariant cyclic homology. 
For later reference we also discuss the notion of an essential module over an idempotented algebra. 
We remark that anti-Yetter-Drinfeld modules are called covariant modules in \cite{Voigtepch}, \cite{Voigtbs}. 
The terminology used here was originally introduced in \cite{HKRS} in the context of Hopf algebras. 
In section \ref{secprim} we discuss the concept of a primitive module over an idempotented algebra and exhibit the 
relation between inductive systems of primitive modules and arbitrary essential modules. 
This is needed for the definition of the local derived category given in section \ref{secpara}. From the point of view of 
homological algebra the local derived category is the main ingredient in the construction of local cyclic homology. 
In section \ref{secantensor} we recall the definition of the analytic tensor algebra and related material from \cite{Meyerthesis}. 
Moreover we review properties of the spectral radius for bornological algebras and discuss locally multiplicative algebras \cite{Meyerborntop}. 
Section \ref{secanalytic} contains the definition 
of the equivariant $ X $-complex of a $ G $-algebra and the definition of equivariant analytic and local cyclic 
homology. This generalizes the constructions in \cite{Meyerthesis}, \cite{Meyerlocalcyclic} as well as the definition 
of entire cyclic cohomology for finite groups given by Klimek, Kondracki and Lesniewski \cite{KKL1}. 
We also discuss briefly the connection to the original approach to local cyclic homology due to Puschnigg. 
In section \ref{sechomstabex} we prove homotopy invariance, stability and excision for 
equivariant analytic and local cyclic homology. The arguments for the analytic and the local theory are analogous
since both theories are constructed in a similar way. In section \ref{seccomp} we study a special situation where 
analytic and local cyclic homology are in fact isomorphic. 
Section \ref{seciso} is devoted to the proof of the isoradial subalgebra theorem. As in the non-equivariant case this theorem is the 
key to establish some nice features of the local theory. In particular, using the isoradial subalgebra 
theorem we study in section \ref{secsmoothing} how local cyclic homology behaves with respect to continuous homotopies and 
stability in the sense of $ C^* $-algebras. 
As a preparation for the definition of the Chern-Connes character in the odd case 
we consider in section \ref{secextprod} the equivariant $ X $-complex of tensor products. 
In section \ref{seckkg} we recall the general approach to bivariant $ K $-theories developped by Cuntz \cite{Cuntzdocumenta}, \cite{Cuntzbivkcstar}. 
Based on the resulting picture of equivariant $ KK $-theory we define the equivariant Chern-Connes character in section \ref{secchernconnes}. 
In the even case the existence of this transformation is an immediate consequence of the universal property of equivariant 
$ KK $-theory \cite{Thomsen}, \cite{Meyerkkg}. As in the non-equivariant case 
the equivariant Chern-Connes character is multiplicative with respect to the Kasparov product and the 
composition product, respectively. Finally, we describe an elementary calculation of the Chern-Connes character in the case of 
profinite groups. More detailed computations together with applications will be discussed in a separate paper. \\ 
Throughout the paper $ G $ will be a second countable totally disconnected locally compact group. All bornological vector spaces 
are assumed to be separated and convex. \\
I am indebted to R. Meyer for providing me his preprint \cite{Meyerlocalcyclic} and answering some questions 
related to local cyclic homology. 

\section{Smooth representations and anti-Yetter-Drinfeld modules}\label{sectotdis}

In this section we recall the basic theory of smooth representations of totally disconnected groups 
and the concept of an anti-Yetter-Drinfeld module. Smooth representations of locally compact groups on bornological 
vector spaces were studied by Meyer in \cite{Meyersmoothrep}. The only difference in our discussion here is that 
we allow for representations on possibly incomplete spaces. Apart from smooth representations, anti-Yetter-Drinfeld modules play a 
central role in equivariant cyclic homology. These modules were called covariant modules in \cite{Voigtepch}, \cite{Voigtbs}. \\
Smooth representations and anti-Yetter-Drinfeld modules for totally disconnected groups can be viewed as essential modules over certain 
idempotented algebras in the following sense. 
\begin{definition} 
An algebra $ H $ with the fine bornology is called idempotented if for every small 
subset $ S \subset H $ there exists an idempotent $ e \in H $ such that $ e \cdot x = x = x \cdot e $ for all 
$ x \in S $. 
\end{definition}
In other words, for every finite set $ F $ of elements in $ H $ there exists an idempotent $ e \in H $ 
which acts like a unit on every element of $ F $. We call a separated $ H $-module $ V $ essential if the 
natural map $ H \otimes_H V \rightarrow V $ is a bornological isomorphism. Since $ H $ carries the fine bornology, the 
completion $ V^c $ of an essential $ H $-module is again essential, and our notion is compatible with the 
concept of an essential module over a bornological algebra with approximate identity \cite{Meyersmoothrep}. 
Clearly an idempotented algebra is a bornological algebra with approximate identity. \\
Let us now consider smooth representations. 
A representation of $ G $ on a separated bornological vector space  $ V $ is a group homomorphism $ \pi: G \rightarrow 
\Aut(V) $ where $ \Aut(V) $ denotes the group of bounded linear automorphisms of $ V $. A bounded linear map 
between representations of $ G $ is called equivariant if it commutes with the action of $ G $. We write 
$ \Hom_G(V,W) $ for the space of equivariant bounded linear maps between the representations $ V $ and $ W $. Let 
$ F(G,V) $ be the space of all functions from $ G $ to $ V $. The adjoint of a representation $ \pi $ is the 
bounded linear map $ [\pi]: V \rightarrow F(G,V) $ given by $ [\pi](v)(t) = \pi(t)(v) $. In the sequel 
we write simply $ t \cdot v $ instead of $ \pi(t)(v) $. \\
We write $ \D(G) $ for the space of smooth functions on $ G $ with compact support equipped with the 
fine bornology. Smoothness of a function $ f $ on a totally disconnected group is equivalent to $ f $ being locally constant. 
If $ V $ is a bornological vector space then $ \D(G) \otimes V = \D(G,V) $ is the space of compactly supported smooth 
functions on $ G $ with values in $ V $. The space $ \E(G,V) $ consists of all smooth functions on $ G $ with values in $ V $.
\begin{definition}
Let $ G $ be a totally disconnected group and let $ V $ be a separated (complete) bornological vector space. A representation 
$ \pi $ of $ G $ on $ V $ is smooth if $ [\pi] $ defines a bounded linear map from $ V $ into $ \E(G,V) $. 
A smooth representation is also called a separated (complete) $ G $-module. 
\end{definition}
Let $ V $ be a separated $ G $-module. Then for every small subset $ S \subset V $ the 
pointwise stabilizer $ G_S $ of $ S $ is an open subgroup of $ G $. 
Conversely, if $ \pi $ is a representation of $ G $ on a bornological vector space $ V $ such that $ G_S $ is 
open for every small subset $ S \subset V $ then $ \pi $ is smooth. 
In particular, if $ V $ carries the fine bornology the above definition reduces to the 
usual definition of a smooth representation on a complex vector space. 
Every representation of a discrete group is smooth. 
Note that a representation $ \pi $ of $ G $ on a separated bornological vector space $ V $ determines a 
representation $ \pi^c $ of $ G $ on the completion $ V^c $. If $ V $ is a separated $ G $-module then $ V^c $ becomes a complete $ G $-module 
in this way. \\
As already mentioned in the beginning, smooth representations can be identified with essential modules over a certain idempotented algebra. 
The Hecke algebra of a totally disconnected group $ G $ is the space 
$ \D(G) $ equipped with the convolution product
$$
(f * g)(t) = \int_G f(s) g(s^{-1}t) ds 
$$
where $ ds $ denotes a fixed left Haar measure on $ G $. Since $ G $ is totally disconnected this algebra is idempotented. 
Every separated $ G $-module $ V $ becomes an essential $ \D(G) $-module by integration,
and conversely, every essential $ \D(G) $-module is obtained in this way. This yields a natural isomorphism 
between the category of separated (complete) $ G $-modules and the category of 
separated (complete) essential $ \D(G) $-modules. \\
A separated (complete) $ G $-algebra is a separated (complete) bornological algebra which is 
also a $ G $-module such that the multiplication $ A \otimes A \rightarrow A $ is equivariant. 
For every separated $ G $-algebra $ A $ the (smooth) crossed product $ A \rtimes G $ is 
the space $ \D(G,A) $ with the convolution multiplication 
$$
(f * g)(t) = \int_G f(s) s \cdot g(s^{-1}t) ds. 
$$
Note in particular that the crossed product associated to the trivial action of $ G $ 
on $ \mathbb{C} $ is the Hecke algebra of $ G $. \\
In connection with actions on $ C^* $-algebras we will have to consider representations of $ G $ which are not smooth. 
For an arbitrary representation of $ G $ on a bornological vector space $ V $ the smoothing $ \Smooth_G(V) $ 
is defined by  
$$
\Smooth_G(V) = \{f \in \E(G,V)| f(t) = t \cdot f(e) \; \text{for all} \; t \in G \}
$$
equipped with the subspace bornology and the right regular representation. We will usually simply write $ \Smooth $ 
instead of $ \Smooth_G $ in the sequel. The smoothing $ \Smooth(V) $ is always a smooth representation of $ G $. If $ V $ is complete, then $ \Smooth(V) $ 
is a complete $ G $-module. There is an injective equivariant bounded linear map $ \iota_V: \Smooth(V) \rightarrow V $ given by 
$ \iota_V(f) = f(e) $. 
\begin{prop} \label{propsmoothing}
Let $ G $ be a totally disconnected group and $ \pi $ be a representation of $ G $ on a separated bornological vector space $ V $. 
The equivariant bounded linear map $ \iota_V: \Smooth(V) \rightarrow V $ induces a natural isomorphism 
$$
\Hom_G(W,V) \cong \Hom_G(W, \Smooth(V))  
$$
for all separated $ G $-modules $ W $. 
\end{prop} 
Hence the smoothing functor $ \Smooth $ is right adjoint to the forgetful functor 
from the category of smooth representations to the category of arbitrary representations. \\
Assume that $ A $ is a separated bornological algebra which is at 
the same time equipped with a representation of $ G $ such that the multiplication $ A \otimes A \rightarrow A $ is equivariant. 
Then $ \Smooth(A) $ is a separated $ G $-algebra in a natural way. This applies in particular to actions on $ C^* $-algebras. 
When $ C^* $-algebras are viewed as bornological algebras we always work with the precompact bornology. 
If $ A $ is a $ G $-$ C^* $-algebra we use the smoothing functor to obtain a complete $ G $-algebra $ \Smooth(A) $. 
We will study properties of this construction in more detail in section \ref{secsmoothing}. \\
Next we discuss the concept of an anti-Yetter-Drinfeld module. Let $ \mathcal{O}_G $ be the 
commutative algebra of compactly supported smooth functions on $ G $ with pointwise multiplication equipped with 
the action of $ G $ by conjugation. 
\begin{definition} Let $ G $ be a totally disconnected group. A separated (complete) $ G $-anti-Yetter-Drinfeld 
module is a separated (complete) bornological vector space $ M $ which is 
both an essential $ \mathcal{O}_G $-module and a $ G $-module such that
\begin{equation*}
s \cdot (f \cdot m) = (s \cdot f) \cdot (s \cdot m)
\end{equation*}
for all $ s \in G, f \in \mathcal{O}_G $ and $ m \in M $. 
\end{definition} 
A morphism $ \phi: M \rightarrow N $ between anti-Yetter-Drinfeld modules is a bounded linear map which is 
$ \mathcal{O}_G $-linear and equivariant. In the sequel we will use the terminology $ \ayd $-module and $ \ayd $-map 
for anti-Yetter-Drinfeld modules and their morphisms. Moreover we denote by $ \SHom_G(M,N) $ the space of $ \ayd $-maps 
between $ \ayd $-modules $ M $ and $ N $. Note that the completion $ M^c $ of a separated $ \ayd $-module $ M $ is a complete $ \ayd $-module. \\
We write $ \AYD(G) $ for the crossed product $ \mathcal{O}_G \rtimes G $. The algebra $ \AYD(G) $ is idempotented and plays the same 
role as the Hecke algebra $ \D(G) $ in the context of smooth representations. More 
precisely, there is an isomorphism of categories between the category of separated (complete) $ \ayd $-modules and 
the category of separated (complete) essential modules over $ \AYD(G) $.  
In particular, $ \AYD(G) $ itself is an $ \ayd $-module in a natural way. We may view elements of $ \AYD(G) $ as 
smooth functions with compact support on $ G \times G $ where the first variable corresponds to $ \mathcal{O}_G $ and 
the second variable corresponds to $ \D(G) $. The multiplication in $ \AYD(G) $ becomes 
$$
(f \cdot g)(s,t) = \int_G f(s,r)g(r^{-1}sr, r^{-1} t) dr 
$$
in this picture. An important feature of this crossed product is that there exists an isomorphism $ T: \AYD(G) \rightarrow \AYD(G) $ of 
$ \AYD(G) $-bimodules given by 
$$
T(f)(s,t) = f(s,st)
$$
for $ f \in \AYD(G) $. More generally, if $ M $ is an arbitrary separated $ \ayd $-module we obtain an 
automorphism of $ M \cong \AYD(G) \otimes_{\AYD(G)} M $ by applying $ T $ to the first tensor factor. By slight abuse 
of language, the resulting map is again denoted by $ T $. This construction is natural in the sense that $ T \phi = \phi T $ for every $ \ayd $-map 
$ \phi: M \rightarrow N $. 

\section{Primitive modules and inductive systems} \label{secprim}

In this section we introduce primitive anti-Yetter-Drinfeld-modules and discuss the relation between inductive systems of primitive modules 
and general anti-Yetter-Drinfeld-modules for totally disconnected groups. This is needed for the definition of equivariant local cyclic homology. \\
Recall from section \ref{sectotdis} that anti-Yetter-Drinfeld modules for a totally disconnected group $ G $ can be viewed as essential modules over 
the idempotented algebra $ \AYD(G) $. Since it creates no difficulties we shall work in the more general setting of essential modules over an 
arbitrary idempotented algebra $ H $ in this section. We let $ \C $ be either the category of separated or complete essential modules over $ H $. Morphisms 
are the bounded $ H $-module maps in both cases. Moreover we let $ \ind(\C) $ be the associated ind-category. The 
objects of $ \ind(\C) $ are inductive systems of objects in $ \C $ and 
the morphisms between $ M = (M_i)_{i \in I} $ and $ (N_j)_{j \in J} $ are given by 
$$
\Hom_{\,\ind(\C)}(M, N) = \varprojlim_{i \in I} \varinjlim_{j \in J} \Hom_\C(M_i, N_j)
$$
where the limits are taken in the category of vector spaces. There is a canonical
functor $ \varinjlim $ from $ \ind(\C) $ to $ \C $ which associates to an inductive system its separated inductive limit. \\
If $ S $ is a small disk in a bornological vector space we write $ \bra S \ket $ for 
the associated normed space. There is a functor which associates to a (complete) bornological vector space $ V $ the 
inductive system of (complete) normed spaces $ \bra S \ket $ where $ S $ runs over the (completant) small disks in $ V $. 
We need a similar construction in the context of $ H $-modules. 
Let $ M $ be a separated (complete) essential $ H $-module and let $ S \subset M $ be a (completant) small disk. 
We write $ H \bra S \ket $ for the image of the natural map $ H \otimes \bra S \ket \rightarrow M $ equipped with 
the quotient bornology and the induced $ H $-module structure. By slight abuse of language we call this module the submodule 
generated by $ S $ and write $ H \bra S \ket \subset M $. 
\begin{definition}
An object of $ \C $ is called primitive if it is generated by a single small disk. 
\end{definition}
In other words, a separated (complete) essential $ H $-module $ P $ is primitive iff there exists a (completant) small disk 
$ S \subset P $ such that the natural map $ H \bra S \ket \rightarrow P $ is an isomorphism. Note that in the special case $ H = \mathbb{C} $ 
the primitive objects are precisely the (complete) normed spaces. \\
Let us write $ \ind(P(\C)) $ for the full subcategory of $ \ind(\C) $ consisting of inductive systems of primitive modules. 
For every $ M \in \C $ we obtain an inductive system of primitive modules over the directed set of 
(completant) small disks in $ M $ by associating to every disk $ S $ the primitive module generated by $ S $. 
This construction yields a functor $ \Analysis $ from $ \C $ to 
$ \ind(P(\C)) $ which will be called the dissection functor. Note that the inductive system $ \Analysis(M) $ 
has injective structure maps for every $ M \in \C $. By definition, an injective inductive system is an inductive system
whose structure maps are all injective. An inductive system is called essentially injective if it is isomorphic 
in $ \ind(\C) $ to an injective inductive system. \\ 
The following assertion is proved in the same way as the corresponding result for bornological vector spaces \cite{Meyerthesis}.  
\begin{prop}\label{lemcovlim}
The direct limit functor $ \varinjlim $ is left adjoint to the dissection functor $ \Analysis $. More precisely, 
there is a natural isomorphism 
$$
\Hom_\C(\varinjlim(M_j)_{j \in J}, N) \cong \Hom_{\,\ind(P(\C))}((M_j)_{j \in J}, \Analysis(N))
$$
for every inductive system $ (M_j)_{j \in J} $ of primitive objects and every $ N \in \C $. 
Moreover 
$ \varinjlim \Analysis $ is naturally equivalent to the identity and the functor $ \Analysis $ is fully faithful. 
\end{prop}
In addition we have that $ \Analysis \varinjlim (M_i)_{i \in I} $ is isomorphic to $ (M_i)_{i \in I} $ provided the system $ (M_i)_{i \in I} $ 
is essentially injective. It follows that the dissection functor $ \Analysis $ induces an equivalence between $ \C $ and the 
full subcategory of $ \ind(P(\C)) $ consisting of all injective inductive systems of primitive modules. 

\section{Paracomplexes and the local derived category} \label{secpara} 

In this section we review the notion of a paracomplex and discuss some related constructions in homological algebra. 
In particular, in the setting of anti-Yetter-Drinfeld modules over a totally disconnected group we define
locally contractible paracomplexes and introduce the local derived category. \\
Let us begin with the definition of a para-additive category \cite{Voigtepch}. 
\begin{definition} A para-additive category is an additive category $ \mathcal{C} $ 
together with a natural isomorphism $ T $ of the identity functor $ \id: \mathcal{C} \rightarrow \mathcal{C} $.   
\end{definition}
It is explained in section \ref{sectotdis} that every $ \ayd $-module is equipped with a natural automorphism denoted by $ T $. Together 
with these automorphisms the category of $ \ayd $-modules becomes a para-additive category in a natural way. In fact, for our purposes this is the main 
example of a para-additive category. 
\begin{definition} Let $ \mathcal{C} $ be a para-additive category. 
A paracomplex $ C = C_0 \oplus C_1 $ in $ \mathcal{C} $ is a given by objects $ C_0 $ and $ C_1 $ 
together with  morphisms $ \partial_0: C_0 \rightarrow C_1 $ and $ \partial_1: C_1 \rightarrow C_0 $ such that
$$ 
\partial^2 = \id - T.
$$
A chain map $ \phi: C \rightarrow D $ between two paracomplexes is a morphism from $ C $ to $ D $ that commutes with
the differentials.
\end{definition}
The morphism $ \partial $ in a paracomplex is called a differential although 
this contradicts the classical definition of a differential. We point out that it does not make sense to speak about the homology of a
paracomplex in general. \\
However, one can define homotopies, mapping cones and suspensions as usual. Moreover, due to 
naturality of $ T $, the space $ \Hom_\C(P,Q) $ of all morphisms between paracomplexes $ P $ and $ Q $ with the standard differential 
is an ordinary chain complex. We write $ {\bf H}(\mathcal{C}) $ for the homotopy category of paracomplexes associated to a 
para-additive category $ \mathcal{C} $. The morphisms in $ {\bf H}(\mathcal{C}) $ are homotopy classes of chain maps. The supension of paracomplexes yields
a translation functor on $ {\bf H}(\mathcal{C}) $. By definition, a triangle 
$$
\xymatrix{
C \; \ar@{->}[r] & X \ar@{->}[r] & Y  \ar@{->}[r] & C[1] \\
 }
$$
in $ {\bf H}(\mathcal{C}) $ is called distinguished if it is isomorphic to a mapping cone triangle. 
As for ordinary chain complexes one proves the following fact. 
\begin{prop}
Let $ \mathcal{C} $ be a para-additive category. Then the homotopy category of 
paracomplexes $ {\bf H}(\mathcal{C}) $ is triangulated. 
\end{prop}
Let us now specialize to the case where $ \mathcal{C} $ is the category of separated (complete) $ \ayd $-modules. 
Hence in the sequel $ {\bf H}(\C) $ will denote the homotopy category of 
paracomplexes of $ \ayd $-modules. We may also consider the homotopy category associated to the corresponding ind-category 
of paracomplexes. There is a direct limit functor $ \varinjlim $ and a 
dissection functor $ \Analysis $ between these categories having the same properties as the 
corresponding functors for $ \ayd $-modules. \\
A paracomplex $ P $ of separated (complete) $ \ayd $-modules is called primitive 
if its underlying $ \ayd $-module is primitive. 
By slight abuse of language, if $ P $ is a primitive paracomplex and $ \iota: P \rightarrow C $ is an injective 
chain map of paracomplexes we will also 
write $ P $ for the image $ \iota(P) \subset C $ with the bornology induced from $ P $. Moreover we 
call $ P \subset C $ a primitive subparacomplex of $ C $ in this case. 
\begin{definition} A paracomplex $ C $ is called locally contractible 
if for every primitive subparacomplex $ P $ of $ C $ the inclusion map $ \iota: P \rightarrow C $ is homotopic to zero. 
A chain map $ f: C \rightarrow D $ between paracomplexes is called a local homotopy equivalence if 
its mapping cone $ C_f $ is locally contractible. 
\end{definition}
The class of locally contractible paracomplexes forms a null system in $ {\bf H}(\C) $. 
We have the following characterization of locally contractible paracomplexes. 
\begin{lemma} \label{primproj}
A paracomplex $ C $ is locally contractible iff $ H_*(\Hom_\C(P,C)) = 0 $ for every 
primitive paracomplex $ P $. 
\end{lemma}
\proof Let $ P \subset C $ be a primitive subparacomplex. If $ H_*(\Hom_\C(P,C)) = 0 $ then 
the inclusion map $ \iota: P \rightarrow C $ is homotopic to zero. It follows that $ C $ 
is locally contractible. Conversely, assume that $ C $ is locally contractible. If $ P $ is a primitive paracomplex 
and $ f: P \rightarrow C $ is a chain map let $ f(P) \subset C $ be the primitive 
subparacomplex corresponding to the image of $ f $. Since $ C $ is locally contractible 
the inclusion map $ f(P) \rightarrow C $ is homotopic to zero. Hence 
the same is true for $ f $ and we deduce $ H_0(\Hom_\C(P,C)) = 0 $. Similarly one obtains 
$ H_1(\Hom_\C(P,C)) = 0 $ since suspensions of primitive paracomplexes are primitive. \qed \\
We shall next construct projective resolutions with respect to the class of locally projective 
paracomplexes. Let us introduce the following terminology. 
\begin{definition} A paracomplex $ P $ is locally projective if $ H_*(\Hom_\C(P, C)) = 0 $ 
for all locally contractible paracomplexes $ C $.
\end{definition}
All primitive paracomplexes are locally projective according to lemma \ref{primproj}. Observe moreover 
that the class of locally projective paracomplexes is closed under direct sums. \\ 
By definition, a locally projective resolution of $ C \in {\bf H}(\C) $ is a locally projective 
paracomplex $ P $ together with a local homotopy equivalence $ P \rightarrow C $. We say that a 
functor $ P: {\bf H}(\C) \rightarrow {\bf H}(\C) $ together with 
a natural transformation $ \pi: P \rightarrow \id $ is a projective resolution functor 
if $ \pi(C): P(C) \rightarrow C $ is a locally projective resolution for all $ C \in {\bf H}(\C) $. 
In order to construct such a functor we proceed as follows. \\
Let $ I $ be a directed set. We view $ I $ as a category with objects the elements of $ I $ 
and morphisms the relations $ i \leq j $. More precisely, there is a morphism $ i \rightarrow j $ from $ i $ to $ j $ in 
this category iff $ i \leq j $. Now consider a functor $ F: I \rightarrow \C $. Such a 
functor is also called an $ I $-diagram in $ \C $. We define a new diagram $ L(F): I \rightarrow \C $ as follows. Set  
\begin{equation*}
L(F)(j) = \bigoplus_{i \rightarrow j} F(i)
\end{equation*}
where the sum runs over all morphisms $ i \rightarrow j $ in $ I $. The map 
$ L(F)(k) \rightarrow L(F)(l) $ induced by a morphism $ k \rightarrow l $ 
sends the summand over $ i \rightarrow k $ identically to the summand over
$ i \rightarrow l $ in $ L(F)(l) $. We have a natural transformation 
$ \pi(F): L(F) \rightarrow F $ sending the summand $ F(i) $ over 
$ i \rightarrow j $ to $ F(j) $ using the map $ F(i \rightarrow j) $. The identical inclusion of 
the summand $ F(j) $ over the identity $ j \rightarrow j $ defines a section $ \sigma(F) $ for $ \pi(F) $. 
Remark that this section is not a natural transformation of $ I $-diagrams in general. \\
Now let $ H: I \rightarrow \C $ be another diagram and let $ (\phi(i): F(i) \rightarrow H(i))_{i \in I} $ be 
an arbitrary family of $ \ayd $-maps. Then there exists a unique natural transformation of $ I $-diagrams
$ \psi: L(F) \rightarrow H $ such that $ \phi(j) = \psi(j) \sigma(F)(j) $ for all $ j $. 
Namely, the summand $ F(i) $ over $ i \rightarrow j $ in $ L(F)(j) $ is mapped under $ \psi(j) $ to $ H(j) $ 
by the map $ H(i \rightarrow j) \phi(i) $. We can rephrase this property as follows. 
Consider the inclusion $ I^{(0)} \subset I $ of all identity morphisms in the category $ I $. 
There is a forgetful functor from the category of $ I $-diagrams to the 
category of $ I^{(0)} $-diagrams in $ \C $ induced by the inclusion $ I^{(0)} \rightarrow I $ and 
a natural isomorphism
$$
\Hom_{I}(L(F), H) \cong \Hom_{I^{(0)}}(F, H)
$$
where $ \Hom_I $ and $ \Hom_{I^{(0)}} $ denote the morphism sets in the categories of 
$ I $-diagrams and $ I^{(0)} $-diagrams, respectively. This means that the previous construction defines a left adjoint 
functor $ L $ to the natural forgetful functor. \\
For every $ j \in I $ we have a split extension of $ \ayd $-modules 
$$
\xymatrix{
   J(F)(j)\; \ar@{>->}[r]^{\iota(F)(j)} & L(F)(j) \ar@{->>}[r]^{\;\;\;\,\pi(F)(j)} & F(j) \\
 }
$$
where by definition $ J(F)(j) $ is the kernel of the $ \ayd $-map $ \pi(F)(j) $ and $ \iota(F)(j) $ 
is the inclusion. 
The $ \ayd $-modules $ J(F)(j) $ assemble to an $ I $-diagram and we obtain an extension 
$$
\xymatrix{
   J(F)\; \ar@{>->}[r]^{\iota(F)} & L(F) \ar@{->>}[r]^{\;\pi(F)} & F \\
 }
$$
of $ I $-diagrams which splits as an extension of $ I^{(0)} $-diagrams. 
We apply the functor $ L $ to the diagram $ J(F) $ and obtain a diagram denoted by $ LJ(F) $ and a 
corresponding extension as before. Iterating this procedure 
yields a family of diagrams $ LJ^n(F) $. 
More precisely, we obtain extensions 
$$
\xymatrix{
   J^{n + 1}(F)\; \ar@{>->}[r]^{\iota(J^n(F))} & LJ^n(F) \ar@{->>}[r]^{\;\;\,\pi(J^n(F))} & J^n(F) \\
 }
$$
for all $ n \geq 0 $ where $ J^0(F) = F $, $ J^1(F) = J(F) $ and $ LJ^0(F) = L(F) $. 
In addition we set $ LJ^{-1}(F) = F $ and $ \iota(J^{-1}(F)) = \id $. By construction there are natural 
transformations $ LJ^n(F) \rightarrow LJ^{n - 1}(F) $ for all $ n $ given by 
$ \iota(J^{n - 1}(F))\pi(J^n(F)) $. In this way we obtain a complex 
$$
\cdots \rightarrow  LJ^3(F) \rightarrow  LJ^2(F) \rightarrow LJ^1(F) \rightarrow LJ^0(F) \rightarrow F \rightarrow 0 
$$
of $ I $-diagrams. Moreover, this complex is split exact as a complex of $ I^{(0)} $-diagrams, that is, 
$ LJ^\bullet(F)(j) $  is a split exact complex of $ \ayd $-modules for all $ j \in I $. \\
Assume now that $ F $ is an $ I $-diagram of paracomplexes in 
$ \C $. We view $ F $ as a pair of $ I $-diagrams $ F_0 $ and $ F_1 $ of $ \ayd $-modules 
together with natural transformations $ \partial_0: F_0 \rightarrow F_1 $ and 
$ \partial_1: F_1 \rightarrow F_0 $ such that $ \partial^2 = \id - T $. Let us construct a family of  $ I $-diagrams
$ (LJ(F),d^h,d^v) $ as follows. Using the same notation as above we set  
$$
LJ(F)_{pq} = LJ^q(F_p) 
$$
for $ q \geq 0 $ and define the horizontal differential $ d^h_{pq}: LJ(F)_{pq} \rightarrow  LJ(F)_{p - 1,q} $ by 
$$ 
d^h_{pq} = (-1)^q LJ^q(\partial_p). 
$$ 
The vertical differential $ d^v_{pq}: LJ(F)_{pq} \rightarrow  LJ(F)_{p,q - 1} $ is 
given by 
$$ 
d^v_{pq} = \iota(J^{q - 1}(F_p)) \pi(J^q(F_p)). 
$$
Then the relations $ (d^v)^2 = 0 $, $ (d^h)^2 = \id - T $ as well as $ d^v d^h + d^h d^v = 0 $ hold. 
Hence, if we define $ \Tot(LJ(F)) $ by 
\begin{equation*}
(\Tot LJ(F))_n = \bigoplus_{p + q = n} LJ(F)_{pq}
\end{equation*}
and equip it with the boundary $ d^h + d^v $ we obtain an $ I $-diagram of paracomplexes. 
We write $ \holim(F) $ for the inductive limit of the diagram $ \Tot LJ(F) $ and call this paracomplex 
the homotopy colimit of the diagram $ F $. There is a canonical
chain map $ \holim(F) \rightarrow \varinjlim(F) $ and a natural filtration on 
$ \holim(F) $ given by 
$$
\holim(F)^{\leq k}_n  = \bigoplus_{\begin{smallmatrix} p + q = n \\ 
q \leq k
\end{smallmatrix}} \varinjlim LJ(F)_{pq}
$$
for $ k \geq 0 $. Observe that the natural inclusion $ \iota^k: \holim(F)^{\leq k} \rightarrow \holim(F) $ is a
chain map and that there is an obvious retraction $ \pi^k: \holim(F) \rightarrow \holim(F)^{\leq k} $ 
for $ \iota^k $. However, this retraction is not a chain map. 
\begin{prop} \label{holimlocproj} 
Let $ F = (F_i)_{i \in I} $ be a directed system of paracomplexes. If the paracomplexes $ F_j $ are 
locally projective then the homotopy colimit $ \emph{ho}\text{-}\varinjlim(F) $ is locally projective as well. 
If the system $ (F_i)_{i \in I} $ is essentially injective then $ \emph{ho}\text{-}\varinjlim(F) \rightarrow \varinjlim(F) $ is
a local homotopy equivalence. 
\end{prop}
\proof Assume first that the paracomplexes $ F_j $ are locally projective. In order to prove that $ \holim(F) $ is locally projective 
let $ \phi: \holim(F) \rightarrow C $ be a chain map 
where $ C $ is a locally contractible paracomplex. We have to show that $ \phi $ is homotopic to zero. The composition 
of the natural map $ \iota^0: \holim(F)^{\leq 0} \rightarrow \holim(F) $ with $ \phi $ 
yields a chain map $ \psi^0 = \phi \iota^0: \varinjlim LJ^0(F) \rightarrow C $. 
By construction of $ LJ^0(F) $ we have isomorphisms
$$
\Hom_\C(\varinjlim LJ^0(F), C) \cong \Hom_I(LJ^0(F),C) \cong \Hom_{I^{(0)}}(F, C)
$$
where we use the notation introduced above and $ C $ is viewed as a constant diagram of paracomplexes. 
Hence, since the paracomplexes $ F_j $ are locally projective, there 
exists a morphism $ h^0 $ of degree one such that $ \partial h^0 + h^0 \partial = \psi^0 $. 
This yields a chain homotopy between $ \psi^0 $ and $ 0 $. 
Using the retraction $ \pi^0: \holim(F) \rightarrow \holim(F)^{\leq 0} $ we obtain 
a chain map $ \phi^1 = \phi - [\partial, h^0 \pi^0] $ from $ \holim(F) $ to $ C $. 
This map is clearly homotopic to $ \phi $ and by construction we have $ \phi^1 \iota^0 = 0 $. Consider next the 
map $ \psi^1 $ given by the composition
$$
\varinjlim LJ^1(F) \rightarrow \holim(F) \rightarrow C 
$$
where the first arrow is the natural one and the second map is $ \phi^1 $. 
Since $ \phi^1 $ vanishes on $ \holim(F)^{\leq 0} $ we see that $ \psi^1 $ 
is a chain map. Observe moreover that $ J^1(F) $ is a locally projective paracomplex. 
The same argument as before yields a homotopy $ h^1: \varinjlim LJ^1(F) \rightarrow C $ such that $ \partial h^1 + h^1 \partial = \psi^1 $. 
We define a chain map $ \phi^2 $ by $ \phi^2 = \phi^1 - [\partial, h^1 \pi^1] = \phi - [\partial, h^1 \pi^1 + h^0 \pi^0] $ and get 
$ \phi^2 \iota^1 = 0 $. Continuing this process we obtain a family of 
$ \ayd $-maps $ h^n: \varinjlim LJ^n(F) \rightarrow C $ which assembles to a homotopy between $ \phi $ and zero. \\
Let $ C(\pi) $ be the mapping cone of the natural map $ \pi: \holim(F) \rightarrow \varinjlim(F) $. 
Moreover we write $ C(j) $ for the mapping cone of $ \Tot LJ(F)(j) \rightarrow F(j) $. 
It follows immediately from the constructions that $ C(j) $ is contractible for every $ j \in I $. Now let $ P \subset C(\pi) $ be a 
primitive subparacomplex. If the system $ F $ is essentially injective then there exists an index 
$ j \in I $ such that $ P \subset C(j) $. Consequently, the map $ P \rightarrow C(\pi) $ is 
homotopic to zero in this case, and we conclude that $ \pi $ is a local homotopy equivalence. \qed \\
Using the previous proposition we can construct a projective resolution functor with respect to 
the class of locally projective paracomplexes. More precisely, one obtains a functor $ P: {\bf H}(\C) \rightarrow {\bf H}(\C) $ by 
setting
$$
P(C) = \holim \Analysis(C) 
$$
for every paracomplex $ C $. In addition, there is a natural transformation $ P \rightarrow \id $ induced by the canonical chain map 
$ \holim(F) \rightarrow \varinjlim(F) $ for every inductive system $ F $. Since $ \Analysis(C) $ is an injective inductive system 
of locally projective paracomplexes for $ C \in {\bf H}(\C) $ 
it follows from proposition \ref{holimlocproj} that this yields a projective resolution functor as desired. \\
Let us now define the local derived category of paracomplexes. 
\begin{definition}
The local derived category $ {\bf D}(\C) $ is the localization of $ {\bf H}(\C) $ with respect to the 
class of locally contractible paracomplexes.
\end{definition}
By construction, there is a canonical functor $ {\bf H}(\C) \rightarrow {\bf D}(\C) $ which sends 
local homotopy equivalences to isomorphisms. Using the projective resolution functor $ P $ one
can describe the morphism sets in the derived category by
$$
\Hom_{{\bf D}(\C)}(C,D) \cong \Hom_{{\bf H}(\C)}(P(C), D) \cong \Hom_{{\bf H}(\C)}(P(C), P(D))
$$
for all paracomplexes $ C $ and $ D $. \\
For the purposes of local cyclic homology we consider the left derived functor of the 
completion functor. This functor is called the derived completion and is given by 
$$
X^{\mathbb{L}c} = P(X)^c 
$$
for every paracomplex $ X $ of separated $ \ayd $-modules. 
Inspecting the construction of the homotopy colimit shows that 
$ X^{\mathbb{L}c} \cong \holim(\Analysis(X)^c) $ where the completion of an inductive system is 
defined entrywise. 

\section{The analytic tensor algebra and the spectral radius} \label{secantensor} 

In this section we discuss the definition of the analytic tensor algebra as well as analytically nilpotent algebras and 
locally multiplicative algebras. The spectral radius of a small subset in a bornological algebra is defined and some of 
its basic properties are established. We refer to \cite{Meyerthesis}, \cite{Meyerlocalcyclic}, \cite{Meyerborntop} for more details. \\
Let $ G $ be a totally disconnected group and let $ A $ be 
a separated $ G $-algebra. We write $ \Omega^n(A) $ for the space of (uncompleted) noncommutative 
$ n $-forms over $ A $. As a bornological vector space one has $ \Omega^0(A) = A $ and
$$ 
\Omega^n(A) = A^+ \otimes A^{\otimes n} 
$$ 
for $ n > 0 $ where $ A^+ $ denotes the unitarization of $ A $. Simple tensors in $ \Omega^n(A) $ 
are usually written in the form $ a_0 da_1 \cdots da_n $ where $ a_0 \in A^+ $ and $ a_j \in A $ for $ j > 0 $. 
Clearly $ \Omega^n(A) $ is a separated $ G $-module with the diagonal action. 
We denote by $ \Omega(A) $ the direct sum of the spaces $ \Omega^n(A) $. The differential $ d $ on 
$ \Omega(A) $ and the multiplication of forms are defined in an obvious way such that the graded Leibniz rule holds. \\
For the purpose of analytic and local cyclic homology it is crucial to consider a bornology on $ \Omega(A) $ 
which is coarser than the standard bornology for a direct sum. 
By definition, the analytic bornology on $ \Omega(A) $ is the bornology generated by the sets 
\begin{equation*}
[S](dS)^\infty = S \cup \bigcup_{n = 1}^\infty S(dS)^n \cup (dS)^n
\end{equation*}
where $ S \subset A $ is small. Here and in the sequel the notation $ [S] $ is used to denote the union of the subset $ S \subset A $ 
with the unit element $ 1 \in A^+ $. Equipped with this bornology $ \Omega(A) $ is again a separated 
$ G $-module. Moreover the differential $ d $ and the multiplication of forms are bounded with respect to the analytic bornology. 
It follows that the Fedosov product defined by 
$$
\omega \circ \eta = \omega \eta - (-1)^{|\omega|} d\omega d\eta
$$ 
for homogenous forms $ \omega $ and $ \eta $ is bounded as well. By definition, the analytic tensor 
algebra $ \mathcal{T}A $ of $ A $ is the 
even part of $ \Omega(A) $ equipped with the Fedosov product and the analytic bornology. 
It is a separated $ G $-algebra in a natural way. Unless explicitly stated otherwise, we will 
always equip $ \Omega(A) $ and $ \mathcal{T}A $ with the analytic bornology in the sequel. \\
The underlying abstract algebra of $ \mathcal{T}A $ can be identified with 
the tensor algebra of $ A $. This relationship between tensor algebras and differential forms is a central 
idea in the approach to cyclic homology developped by Cuntz and Quillen \cite{CQ1}, \cite{CQ2}, \cite{CQ4}. 
However, since the analytic bornology is different from the direct sum bornology, the analytic 
tensor algebra $ \mathcal{T}A $ is no longer universal for all equivariant bounded linear maps from $ A $ 
into separated $ G $-algebras. In order to formulate its universal property correctly 
we need some more terminology. \\
The curvature of an equivariant bounded linear map $ f: A \rightarrow B $ between separated $ G $-algebras is the 
equivariant linear map $ \omega_f: A \otimes A \rightarrow B $ 
given by 
$$
\omega_f(x,y) = f(xy) - f(x) f(y). 
$$
By definition, the map $ f $ has analytically nilpotent curvature if 
$$
\omega_f(S,S)^\infty = \bigcup_{n = 1}^\infty \omega_f(S,S)^n 
$$
is a small subset of $ B $ for all small subsets $ S \subset A $. An equivariant bounded linear map $ f: A \rightarrow B $ with analytically nilpotent 
curvature is called an equivariant lanilcur. 
The analytic bornology is defined in such a way that the equivariant homomorphism $ [[f]]: \mathcal{T}A \rightarrow B $ 
associated to an equivariant bounded linear map $ f: A \rightarrow B $ is 
bounded iff $ f $ is a lanilcur. \\
It is clear that every bounded homomorphism $ f: A \rightarrow B $ is a lanilcur. In particular, the identity map of $ A $ corresponds to the bounded homomorphism 
$ \tau_A: \mathcal{T}A \rightarrow A $ given by the canonical projection onto 
differential forms of degree zero. The kernel of the map $ \tau_A $ is denoted by $ \mathcal{J}A $, and we obtain an extension 
$$
\xymatrix{
   \mathcal{J}A\; \ar@{>->}[r] & \mathcal{T}A \ar@{->>}[r] & A \\
}
$$
of separated $ G $-algebras. This extension has an equivariant bounded linear splitting $ \sigma_A $ given by 
the inclusion of $ A $ as differential forms of degree zero. The algebras $ \mathcal{J}A $ and $ \mathcal{T}A $ 
have important properties that we shall discuss next. \\
A separated $ G $-algebra $ N $ is called analytically nilpotent if 
$$
S^\infty = \bigcup_{n \in \mathbb{N}} S^n
$$
is small for all small subsets $ S \subset N $. For instance, every nilpotent bornological algebra is 
analytically nilpotent. The ideal $ \mathcal{J}A $ in the analytic tensor algebra of a bornological $ A $ is an important example of 
an analytically nilpotent algebra. \\
A separated $ G $-algebra $ R $ is called equivariantly analytically quasifree provided the following condition is 
satisfied. If $ K $ is an analytically nilpotent $ G $-algebra and 
$$
\xymatrix{
   K \; \ar@{>->}[r] & E \ar@{->>}[r] & Q \\
}
$$
is an extension of complete $ G $-algebras with equivariant bounded linear splitting then for 
every bounded equivariant homomorphism $ f: R \rightarrow Q $ there exists a bounded equivariant lifting 
homomorphism $ F: R \rightarrow E $. The analytic tensor algebra $ \mathcal{T}A $ of a $ G $-algebra $ A $ is a basic example of an 
equivariantly analytically quasifree $ G $-algebra. Another fundamental example is given by the algebra $ \mathbb{C} $ with the 
trivial action. Every equivariantly analytically quasifree $ G $-algebra is in particular equivariantly quasifree in the sense of \cite{Voigtepch}. \\ 
We shall next discuss the concept of a locally multiplicative $ G $-algebra. If $ A $ is a bornological algebra then a disk $ T \subset A $ is called 
multiplicatively closed provided $ T \cdot T \subset T $. 
A separated bornological algebra $ A $ is called locally multiplicative if for every small subset $ S \subset A $ there 
exists a positive real number $ \lambda $ and a small multiplicatively closed disk $ T \subset A $ such that $ S \subset \lambda T $.
It is easy to show that a separated (complete) bornological algebra is locally multiplicative iff it is a direct limit of 
(complete) normed algebras. We point out that the group action on a $ G $-algebra usually does not 
leave multiplicatively closed disks invariant. In particular, a locally multiplicative $ G $-algebra can not be written as a direct limit 
of normed $ G $-algebras in general. \\
It is clear from the definitions that analytically nilpotent algebras are locally multiplicative. In fact, 
locally multiplicatively algebras and analytically nilpotent algebras can be characterized in a concise way 
using the notion of spectral radius. 
\begin{definition}
Let $ A $ be a separated bornological algebra and let $ S \subset A $ be a small subset. The 
spectral radius $ \rho(S) = \rho(S;A) $ is the infimum of all positive real numbers 
$ r $ such that 
\begin{equation*}
(r^{-1} S)^\infty = \bigcup_{n = 1}^\infty (r^{-1} S)^n 
\end{equation*}
is small. If no such number $ r $ exists set $ \rho(S) = \infty $.
\end{definition}
A bornological algebra $ A $ is locally multiplicative 
iff $ \rho(S) < \infty $ for all small subsets $ S \subset A $. Similarly, a bornological algebra is analytically 
nilpotent iff $ \rho(S) = 0 $ for all small subsets $ S $ of $ A $. \\ 
Let us collect some elementary properties of the spectral radius. 
If $ \lambda $ is a positive real number then $ \rho(\lambda S) = \lambda \rho(S) $ for every small subset $ S $. 
Moreover one has $ \rho(S^n) = \rho(S)^n $ for all $ n > 0 $. Remark also that the spectral radius does not distinguish between a 
small set and its disked hull. Finally, let $ f: A \rightarrow B $ be a bounded homomorphism 
and let $ S \subset A $ be small. Then the spectral radius is contractive in the sense that 
$$ 
\rho(f(S);B) \leq \rho(S;A) 
$$ 
since $ f((r^{-1} S)^\infty) = (r^{-1}f(S))^\infty \subset B $ is small 
provided $ (r^{-1} S)^\infty $ is small. 

\section{Equivariant analytic and local cyclic homology} \label{secanalytic}

In this section we recall the definition of equivariant differential forms and the equivariant $ X $-complex and define 
equivariant analytic and local cyclic homology. In addition we discuss the relation to equivariant entire cyclic homology for 
finite groups in the sense of Klimek, Kondracki and Lesniewski and the original definition of local cyclic homology due to Puschnigg. \\
First we review basic properties of equivariant differential forms. The equivariant $ n $-forms over a separated $ G $-algebra $ A $ 
are defined by $ \Omega^n_G(A) = \mathcal{O}_G \otimes \Omega^n(A) $ where $ \Omega^n(A) $ 
is the space of uncompleted differential $ n $-forms over $ A $. The group $ G $ acts diagonally on $ \Omega^n_G(A) $ and we have 
an obvious $ \mathcal{O}_G $-module structure given by multiplication on the first tensor factor. 
In this way the space $ \Omega^n_G(A) $ becomes a separated $ \ayd $-module. \\
On equivariant differential forms we consider the following operators. We have the differential $ d: \Omega^n_G(A) \rightarrow \Omega^{n + 1}_G(A) $ given by 
$$
d(f(s) \otimes x_0 dx_1 \cdots dx_n) = f(s) \otimes dx_0 dx_1 \cdots dx_n 
$$
and the equivariant Hochschild boundary $ b: \Omega^n_G(A) \rightarrow \Omega^{n - 1}_G(A) $ defined by 
\begin{align*}
b(f(s) \otimes &x_0 dx_1 \cdots dx_n) = f(s) \otimes x_0 x_1 dx_2 \cdots dx_n \\
&+ \sum_{j = 1}^{n - 1} (-1)^j f(s) \otimes x_0 dx_1 \cdots d(x_j x_{j + 1}) \cdots dx_n \\
&+ (-1)^n f(s) \otimes (s^{-1} \cdot x_n)x_0 dx_1 \cdots dx_{n - 1}.
\end{align*}
Moreover there is the equivariant Karoubi operator $ \kappa: \Omega^n_G(A) \rightarrow \Omega^n_G(A) $ and 
the equivariant Connes operator $ B: \Omega^n_G(A) \rightarrow \Omega^{n + 1}_G(A) $ which are 
given by the formulas
\begin{equation*}
\kappa(f(s) \otimes x_0dx_1 \cdots dx_n) = (-1)^{n - 1}
f(s) \otimes (s^{-1} \cdot dx_n) x_0 dx_1 \cdots dx_{n - 1}
\end{equation*}
and 
\begin{equation*}
B(f(s) \otimes x_0dx_1 \cdots dx_n) = \sum_{i = 0}^n (-1)^{ni}
f(s) \otimes s^{-1} \cdot(dx_{n + 1 - i} \cdots dx_n)dx_0 \cdots dx_{n - i},
\end{equation*}
respectively. All these operators are $ \ayd $-maps, and the natural symmetry operator $ T $ for $ \ayd $-modules is of the form
\begin{equation*}
T(f(s) \otimes \omega) = f(s) \otimes s^{-1} \cdot \omega 
\end{equation*}
on equivariant differential forms. We shall write $ \Omega_G(A) $ for the direct sum of the spaces $ \Omega^n_G(A) $ in the sequel. 
The analytic bornology on $ \Omega_G(A) $ is defined using the identification $ \Omega_G(A) = \mathcal{O}_G \otimes \Omega(A) $. \\
Together with the operators $ b $ and $ B $ the space $ \Omega_G(A) $ of equivariant differential forms may be viewed as a paramixed complex \cite{Voigtepch}
which means that the relations $ b^2 = 0 $, $ B^2 = 0 $ and 
$ [b,B] = bB + Bb = \id - T $ hold. An important purpose for which equivariant differential forms are needed 
is the definition of the equivariant $ X $-complex of a $ G $-algebra. 
\begin{definition} Let $ A $ be a separated $ G $-algebra. The equivariant $ X $-complex 
$ X_G(A) $ of $ A $ is the paracomplex   
\begin{equation*}
    \xymatrix{
      {X_G(A) \colon \ }
        {\Omega^0_G(A)\;} \ar@<1ex>@{->}[r]^-{d} &
          {\;\Omega^1_G(A)/ b(\Omega^2_G(A)).} 
            \ar@<1ex>@{->}[l]^-{b} 
               } 
\end{equation*}
\end{definition}
Remark in particular that if $ \partial $ denotes the boundary operator in $ X_G(A) $ then the relation $ \partial^2 = \id - T $ 
follows from the fact that equivariant differential forms are a paramixed complex. \\
After these preparations we come to the definition of equivariant analytic cyclic homology. 
\begin{definition} \label{defanalytic} 
Let $ G $ be a totally disconnected group and let $ A $ and $ B $ be separated $ G $-algebras. 
The bivariant equivariant analytic cyclic homology of $ A $ and $ B $ is
\begin{equation*}
HA^G_*(A,B) =
H_*(\SHom_G(X_G(\mathcal{T}(A \otimes \mathcal{K}_G))^c, 
X_G(\mathcal{T}(B \otimes \mathcal{K}_G))^c)).
\end{equation*}
\end{definition}
The algebra $ \mathcal{K}_G $ occuring in this definition is the subalgebra of the algebra of compact operators $ \mathbb{K}(L^2(G)) $ on the 
Hilbert space $ L^2(G) $ obtained as the linear span of all rank-one operators $ |\xi \ket \bra \eta| $ with 
$ \xi, \eta \in \D(G) $. This algebra is equipped with the fine bornology and the action induced from $ \mathbb{K}(L^2(G)) $. 
An important property of the $ G $-algebra $ \mathcal{K}_G $ is that it is projective as a $ G $-module. \\
We point out that the $ \Hom $-complex on the right hand side of the definition, equipped with the usual boundary operator, is an ordinary chain 
complex although both entries are only paracomplexes. Remark also that for the trivial group one reobtains the definition of analytic cyclic homology 
given in \cite{Meyerthesis}. \\
It is frequently convenient to replace the paracomplex $ X_G(\mathcal{T}(A \otimes \mathcal{K}_G)) $ in the definition of the 
analytic theory with another paracomplex construced using the standard boundary $ B + b $ in cyclic homology. 
For every separated $ G $-algebra $ A $ there is a natural isomorphism $ X_G(\mathcal{T}A) \cong \Omega^{\an}_G(A) $ 
of $ \ayd $-modules where $ \Omega^{\an}_G(A) $ is the space $ \Omega_G(A) $ equipped with the 
transposed analytic bornology. The transposed analytic bornology is the bornology generated by the sets 
\begin{equation*}
D \otimes S \cup D \otimes [S] dS \cup \bigcup_{n = 1}^\infty n!\, D \otimes [S][dS](dS)^{2n} 
\end{equation*}
where $ D \subset \mathcal{O}_G $ and $ S \subset A $ are small. 
The operators $ b $ and $ B $ are bounded with respect to the transposed analytic bornology. It follows that
$ \Omega^{\an}_G(A) $ becomes a paracomplex with the differential $ B + b $. 
We remark that rescaling with the constants $ n! $ in degree $ 2n $ and $ 2n + 1 $ yields an isomorphism between $ \Omega^{\an}_G(A) $ 
and the space $ \Omega_G(A) $ equipped with the analytic bornology. 
\begin{theorem}\label{homotopyeq}
Let $ G $ be a totally disconnected group. For every separated $ G $-algebra $ A $ there exists a bornological 
homotopy equivalence between the paracomplexes $ X_G(\mathcal{T}A) $ and $ \Omega^{\an}_G(A) $. 
\end{theorem}
\proof The proof follows the one for the equivariant periodic theory \cite{Voigtepch} and the corresponding 
assertion in the nonequivariant situation \cite{Meyerthesis}. Let $ Q_n: \Omega_G(A) \rightarrow \Omega^n_G(A) \subset \Omega_G(A) $ be the canonical 
projection. Using the explicit formula for the Karoubi operator one checks that the set 
$ \{ C^n \kappa^j Q_n | \, 0\leq j \leq n, n \geq 0 \, \} $ of operators is equibounded on $ \Omega_G(A) $ 
with respect to the analytic bornology for every $ C \in \mathbb{R} $. 
Similarly, the set $ \{\kappa^n Q_n |n \geq 0 \} $ is equibounded with respect to the analytic bornology and 
hence $ \{ C^n \kappa^j Q_n | \, 0\leq j \leq kn, n \geq 0 \, \} $ is equibounded as well for each $ k \in \mathbb{N} $. 
Thus an operator on $ \Omega_G(A) $ of the form $ \sum_{j = 0}^\infty Q_n h_n(\kappa) $ is bounded 
with respect to the analytic bornology if $ (h_n)_{n \in \mathbb{N}} $ is a sequence of polynomials 
whose degrees grow at most linearly and whose absolut coefficient sums grow at most exponentially. By definition,
the absolute coefficient sum of $ \sum_{j = 0}^k a_j x^j $ is $ \sum_{j = 0}^k |a_j| $.  
The polynomials $ f_n $ and $ g_n $ occuring in the proof of theorem 8.6 in \cite{Voigtepch} satisfy these conditions. 
Based on this observation, a direct inspection shows that the maps involved in the definition of the desired homotopy equivalence 
in the periodic case induce bounded maps on $ \Omega_G(A) $ with respect to the analytic bornology. This yields the assertion. \qed \\
Let $ G $ be a finite group and let $ A $ be a unital Banach algebra on which $ G $ 
acts by bounded automorphisms. Klimek, Kondracki and Lesniewski defined the equivariant entire cyclic cohomology of $ A $ in 
this situation \cite{KKL1}. We may also view $ A $ as a bornological algebra with the bounded bornology and consider 
the equivariant analytic theory of the resulting $ G $-algebra.
\begin{prop}
Let $ G $ be a finite group acting on a unital Banach algebra $ A $ by bounded automorphisms. Then the equivariant entire cyclic cohomology of $ A $ 
coincides with the equivariant analytic cyclic cohomology $ HA^G_*(A, \mathbb{C}) $ where $ A $ is viewed as a $ G $-algebra 
with the bounded bornology.
\end{prop}
\proof It will be shown in proposition \ref{defcomp} below that tensoring with the algebra $ \mathcal{K}_G $ is not needed in the definition of 
$ HA^G_* $ for finite groups. Let us write $ C(G) $ for the space of functions on the finite group $ G $. Using 
theorem \ref{homotopyeq} we see that the analytic cyclic cohomology $ HA^G_*(A, \mathbb{C}) $ is computed by the complex consisting of 
families $ (\phi_n)_{n \geq 0} $ of $ n + 1 $-linear maps $ \phi_n: A^+ \times A^n \rightarrow C(G) $ which are equivariant in the sense 
that 
$$
\phi_n(t \cdot a_0, t \cdot a_1, \dots, t \cdot a_n)(s) = \phi_n(a_0, a_1, \dots, a_n)(t^{-1} s t)
$$
and satisfy the entire growth condition 
$$
[n/2]! \; \max_{t \in G} |\phi_n(a_0, a_1, \dots, a_n)(t)| \leq c_S 
$$
for $ a_0 \in [S], a_1, \dots, a_n \in S $ and all small sets $ S $ in $ A $. Here $ [n/2] = k $ for $ n = 2k $ or $ n = 2k + 1 $ and $ c_S $ is a constant 
depending on $ S $. The boundary operator is induced by $ B + b $. An argument analogous to the one due to Khalkhali in 
the non-equivariant case \cite{Khalkhali} shows that this complex is homotopy equivalent to the complex used by Klimek, Kondracki and Lesniewski. \qed 
\begin{definition} \label{deflocal}
Let $ G $ be a totally disconnected group and let $ A $ and $ B $ be separated $ G $-algebras. 
The bivariant equivariant local cyclic homology $ HL^G_*(A,B) $ of $ A $ and $ B $ is
given by 
\begin{equation*}
H_*(\SHom_G(X_G(\mathcal{T}(A \otimes \mathcal{K}_G))^{\mathbb{L}c}, X_G(\mathcal{T}(B \otimes \mathcal{K}_G))^{\mathbb{L}c})).
\end{equation*}
\end{definition}
Recall that the derived completion $ X^{\mathbb{L}c} $ of 
a paracomplex $ X $ was introduced in section \ref{secpara}. 
In terms of the local derived category of paracomplexes defininition \ref{deflocal} can be reformulated in the following way. 
The construction of the derived completion shows together with proposition \ref{holimlocproj} that the paracomplex
$ X_G(\mathcal{T}(A \otimes \mathcal{K}_G))^{\mathbb{L}c} $ is locally projective for every separated $ G $-algebra $ A $. 
It follows that the local cyclic homology group $ HL^G_0(A,B) $ is equal to the space of morphisms 
in the local derived category between $ X_G(\mathcal{T}(A \otimes \mathcal{K}_G))^{\mathbb{L}c} $ 
and $ X_G(\mathcal{T}(B \otimes \mathcal{K}_G))^{\mathbb{L}c} $. 
Consequently, the passage from the analytic theory to the local theory consists in passing from the homotopy category of paracomplexes to 
the local derived category and replacing the completion functor by the derived completion. \\
Both equivariant analytic and local cyclic homology are equipped with an obvious composition product. Every bounded equivariant 
homomorphism $ f: A \rightarrow B $ induces an element $ [f] $ in $ HA^G_*(A,B) $ and in $ HL^G_*(A,B) $, respectively. In particular, 
the identity map $ \id: A \rightarrow A $ defines an element in these theories which acts as a unit with respect to the composition product. \\ 
If $ G $ is the trivial group then definition \ref{deflocal} reduces to the local cyclic theory defined by Meyer in 
\cite{Meyerlocalcyclic}. Let us briefly mention how this definition of local cyclic homology is related to the original 
approach by Puschnigg. In \cite{Puschnigg4} a Fr\'echet algebra $ A $ is called nice if there is a neighborhood of the 
origin $ U $ such that $ S^\infty $ is precompact for all compact sets $ S \subset U $. This condition 
is equivalent local multiplicativity if $ A $ is viewed as a bornological algebra with the precompact 
bornology \cite{Meyerthesis}. 
Hence the class of nice Fr\'echet algebras can be viewed as a particular class of locally multiplicative bornological algebras. 
\begin{prop}
Let $ (A_i)_{i \in I} $ and $ (B_j)_{j \in J} $ be inductive systems of nice Fr\'echet algebras and let $ A $ and $ B $
denote their direct limits, respectively. If the systems $ (A_i)_{i \in I} $ and $ (B_j)_{j \in J} $ have injective structure 
maps then $ HL_*(A,B) $ is naturally isomorphic to the bivariant local cyclic homology for $ (A_i)_{i \in I} $ and $ (B_j)_{j \in J} $ 
as defined by Puschnigg.
\end{prop}
\proof According to the assumptions the inductive system $ \Analysis(X(\mathcal{T}A)) $ is isomorphic to the formal 
inductive limit of $ \Analysis(X(\mathcal{T}A_i))_{i \in I} $ in the category of inductive systems of 
complexes. The completion of the latter is equivalent to the inductive 
system that is used in \cite{Puschnigg4} to define the local theory. Comparing the construction of the local derived category with the definition of the 
derived ind-category given by Puschnigg yields the assertion. \qed \\
Consequently, the main difference between their approaches is that Meyer works explicitly in the setting of bornological vector spaces whereas Puschnigg uses 
inductive systems and considers bornologies only implicitly. 
It should be pointed out that in practice one is mainly interested in inductive systems of nice Fr\'echet algebras with
injective structure maps. Remark moreover that the definition of Meyer also applies 
to bornological algebras which are not locally multiplicative. 

\section{Homotopy invariance, stability and excision}\label{sechomstabex}

In this section we show that equivariant analytic and local cyclic homology are invariant under smooth equivariant 
homotopies, stable and satisfy excision in both variables. \\
For the proof of homotopy invariance and stability of the local theory we need some information about partial completions. 
A subset $ \mathcal{V} $ of a bornological vector space $ V $ is called locally dense if for any small subset $ S \subset V $ 
there is a small disk $ T \subset V $ such that any $ v \in S $ is the limit of a $ T $-convergent sequence 
with entries in $ \mathcal{V} \cap T $. If $ V $ is a metrizable locally convex vector space endowed with the precompact bornology then a 
subset $ \mathcal{V} \subset V $ is locally dense iff it is dense in $ V $ in the topological sense \cite{Meyerborntop}. 
Let $ \mathcal{V} $ be a bornological vector space and let $ i: \mathcal{V} \rightarrow V $ be a bounded linear map 
into a separated bornological vector space $ V $. Then $ V $ together with the map $ i $ is called a partial completion of $ \mathcal{V} $ if 
$ i $ is a bornological embedding and has locally dense range. \\
We will need the following property of partial completions. 
\begin{lemma} \label{parcomplem}
Let $ i: \mathcal{A} \rightarrow A $ be a partial completion of separated $ G $-algebras. Then the induced chain
map $ X_G(\mathcal{T}\mathcal{A})^{\mathbb{L}c} \rightarrow X_G(\mathcal{T}A)^{\mathbb{L}c} $ 
is an isomorphism. If the derived completion is replaced by the ordinary completion the corresponding chain
map is an isomorphism as well. 
\end{lemma}
\proof Let us abbreviate $ C = X_G(\mathcal{T}\mathcal{A}) $ and $ D = X_G(\mathcal{T}A) $. 
It suffices to show that the natural map $ \Analysis(C)^c \rightarrow \Analysis(D)^c $ is an isomorphism 
of inductive systems. 
Since $ i: \mathcal{A} \rightarrow A $ is a partial completion the same holds true for the induced 
chain map $ C \rightarrow D $. 
By local density, for any small disk $ S \subset D $ there exists a small disk $ T \subset D $ 
such that any point in $ S $ is the limit of a $ T $-convergent sequence with entries in $ C \cap T $. 
Observe that $ C \cap T $ is a small disk in $ C $ since the inclusion is a bornological embedding. Consider the 
isometry $ \bra C \cap T \ket \rightarrow \bra T \ket $. By construction, the space $ \bra S \ket $ is 
contained in the range of the isometry $ \bra C \cap T \ket^c \rightarrow \bra T \ket^c $ obtained by applying 
the completion functor. Since $ \bra C \cap T \ket^c $ maps naturally into $ (\AYD(G)\bra C \cap T \ket)^c $ we 
get an induced $ \ayd $-map $ \AYD(G) \bra S \ket \rightarrow (\AYD(G) \bra C \cap T \ket)^c $. Using this observation one 
checks easily that the completions of the inductive systems $ \Analysis(C) $ and $ \Analysis(D) $ are isomorphic. \qed \\
We refer to \cite{Meyerborntop} for the definition of smooth functions with values in a bornological vector 
space. For metrizable locally convex vector spaces with the precompact bornology one reobtains the usual notion. 
Let $ B $ be a separated $ G $-algebra and denote by $ C^\infty([0,1], B) $ the $ G $-algebra of smooth functions on the interval 
$ [0,1] $ with values in $ B $. The group $ G $ acts pointwise on functions, and if $ B $ is complete there is a natural 
isomorphism $ C^\infty([0,1], B) \cong C^\infty[0,1] \cotimes B $. 
A smooth equivariant homotopy is a bounded equivariant homomorphism 
$ \Phi: A \rightarrow C^\infty([0,1],B) $. Evaluation at $ t \in [0,1] $ yields an equivariant homomorphism
$ \Phi_t: A \rightarrow B $. Two equivariant homomorphisms from $ A $ to $ B $ are called equivariantly homotopic 
if they can be connected by an equivariant homotopy. 
\begin{prop}[Homotopy invariance] \label{homotopyinv} Let $ A $ and $ B $ be 
separated $ G $-algebras and let $ \Phi: A \rightarrow C^\infty([0,1],B) $ be a smooth equivariant homotopy. Then 
the induced elements $ [\Phi_0] $ and $ [\Phi_1] $ in $ HL^G_*(A,B) $ are equal. 
An analogous statement holds for the analytic theory.
Hence $ HA^G_* $ and $ HL^G_* $ are homotopy invariant in both variables with respect to smooth equivariant homotopies.  
\end{prop}
\proof For notational simplicity we shall suppress occurences of the algebra $ \mathcal{K}_G $ in our notation. 
Assume first that the homotopy $ \Phi $ is a map from 
$ A $ into $ \mathbb{C}[t] \otimes B $ where 
$ \mathbb{C}[t] $ is viewed as a subalgebra of $ C^\infty[0,1] $ with the subspace bornology. 
The map $ \Phi $ 
induces a bounded equivariant homomorphism $ \mathcal{T}A \rightarrow 
\mathbb{C}[t] \otimes \mathcal{T}B $ since the algebra $ C^\infty[0,1] $ is locally multiplicative. 
As in the proof of homotopy invariance for equivariant periodic cyclic homology ~\cite{Voigtepch} we 
see that the chain maps $ X_G(\mathcal{T}A) \rightarrow X_G(\mathcal{T}B) $ 
induced by $ \Phi_0 $ and $ \Phi_1 $ are homotopic. Consider in particular the 
equivariant homotopy $ \Phi: \mathbb{C}[x] \otimes B \rightarrow \mathbb{C}[t] \otimes \mathbb{C}[x] \otimes B $ defined 
by $ \Phi(p(x) \otimes b) = p(tx) \otimes b $. We deduce that the map $ B \rightarrow C[x] \otimes B $ 
that sends $ b $ to $ b \otimes 1 $ induces 
a homotopy equivalence between $ X_G(\mathcal{T}(\mathbb{C}[x] \otimes B)) $ and $ X_G(\mathcal{T}B) $. 
It follows in particular that the chain maps $ X_G(\mathcal{T}(\mathbb{C}[x] \otimes B)) \rightarrow 
X_G(\mathcal{T}B) $ given by evaluation at $ 0 $ and $ 1 $, respectively, are homotopic. \\
Let us show that $ \mathbb{C}[t] \otimes B \rightarrow C^\infty([0,1],B) $ 
is a partial completion. It suffices to consider the corresponding map for a normed subspace $ V \subset B $ since 
source and target of this map are direct limits of the associated inductive systems with injective structure maps. 
For a normed space $ V $ the assertion follows from Grothendieck's description of 
bounded subsets of the projective tensor product $ C^\infty[0,1] \cotimes_\pi V $. \\
Due to lemma \ref{parcomplem} the chain map $ X_G(\mathcal{T}(\mathbb{C}[t] \otimes B))^{\mathbb{L}c} 
\rightarrow X_G(\mathcal{T}C^\infty([0,1], B))^{\mathbb{L}c} $
is an isomorphism. Hence the chain maps
$
X_G(\mathcal{T}C^\infty([0,1],B))^{\mathbb{L}c} \rightarrow X_G(\mathcal{T}B)^{\mathbb{L}c}
$ 
induced by evalution at $ 0 $ and $ 1 $ are homotopic as well. Now let $ \Phi: A \rightarrow C^\infty([0,1], B) $ 
be an arbitrary homotopy. According to our previous argument, composing the induced chain map 
$ X_G(\mathcal{T}A)^{\mathbb{L}c} \rightarrow X_G(\mathcal{T}C^\infty([0,1], B))^{\mathbb{L}c} $ with the evaluation maps at $ 0 $ and $ 1 $ yields 
the claim for the local theory. The assertion for the analytic theory are obtained in the same way. \qed \\
Next we study stability. Let $ V $ and $ W $ be separated $ G $-modules and let 
$ b: W \times V \rightarrow \mathbb{C} $ be an equivariant bounded bilinear map. Then 
$ l(b) = V \otimes W $ is a separated $ G $-algebra with multiplication 
\begin{equation*}
(x_1 \otimes y_1) \cdot (x_2 \otimes y_2) = x_1 \otimes b(y_1, x_2) y_2 
\end{equation*}
and the diagonal $ G $-action. A particular example is the algebra 
$ \mathcal{K}_G $ which is obtained using the left regular representation $ V = W = \D(G) $ and the pairing 
\begin{equation*}
b(f,g) = \int_G f(s) g(s) ds
\end{equation*}
with respect to left Haar measure. \\
Let $ V $ and $ W $ be separated $ G $-modules and let $ b: W \times V $ be an 
equivariant bounded bilinear map. The pairing $ b $ is called admissible if 
there exists nonzero $ G $-invariant vectors $ v \in V $ and $ w \in W $ such that $ b(w,v) = 1 $. 
In this case $ p = v \otimes w $ is an invariant idempotent element in $ l(b) $ and there is an equivariant homomorphism 
$ \iota_A: A \rightarrow A \otimes l(b) $ given by $ \iota_A(a) = a \otimes p $. 
\begin{prop}\label{StabLemma} Let $ A $ be a separated $ G $-algebra and
let $ b: W \times V \rightarrow \mathbb{C} $ be an admissible pairing. Then
the map $ \iota_A $ induces a homotopy equivalence
$ X_G(\mathcal{T}A)^{\mathbb{L}c} \simeq X_G(\mathcal{T}(A \otimes l(b)))^{\mathbb{L}c} $. 
If the derived completion is replaced by the ordinary completion the corresponding map is 
a homotopy equivalence as well. 
\end{prop}
This result is proved in the same way as in \cite{Voigtepch} using homotopy invariance.
As a consequence we obtain the following stability properties of equivariant analytic and 
local cyclic homology.
\begin{prop}[Stability] \label{stability} Let $ A $ be a separated $ G $-algebra and let
$ b: W \times V \rightarrow \mathbb{C} $ be a nonzero equivariant bounded bilinear map. 
Moreover let $ l(b,A) $ be any partial completion of $ A \otimes l(b) $. Then there exist
invertible elements in $ HL^G_0(A, l(b,A)) $ and $ HA^G_0(A, l(b, A)) $.
\end{prop}
\proof For the uncompleted stabilization $ A \otimes l(b) $ the argument for 
the periodic theory in \cite{Voigtepch} carries over. 
If $ l(b,A) $ is a partial completion of $ A \otimes l(b) $
the natural chain map $ X_G(\mathcal{T}(A \otimes l(b) \otimes \mathcal{K}_G)) \rightarrow 
X_G(\mathcal{T}(l(b, A) \otimes \mathcal{K}_G)) $ becomes an isomorphism after applying the 
(left derived) completion functor according to lemma \ref{parcomplem}. \qed \\
An application of theorem \ref{StabLemma} yields a simpler description of $ HA^G_* $ and $ HL^G_* $ in the case that $ G $ 
is a profinite group. If $ G $ is compact the trivial one-dimensional representation is contained in 
$ \D(G) $. Hence the pairing used to define the algebra $ \mathcal{K}_G $ is admissible in this case. This 
implies immediately the following assertion. 
\begin{prop} \label{defcomp} Let $ G $ be a compact group. Then we have 
a natural isomorphism
\begin{equation*}
HL^G_*(A,B) \cong H_*(\SHom_G(X_G(\mathcal{T}A)^{\mathbb{L}c}, X_G(\mathcal{T}B)^{\mathbb{L}c}))
\end{equation*}
for all separated $ G $-algebras $ A $ and $ B $. An analogous statement holds for the analytic theory. 
\end{prop}
To conclude this section we show that equivariant analytic and local cyclic homology satisfy 
excision in both variables. 
\begin{theorem}[Excision]\label{Excision} Let $ A $ be a separated $ G $-algebra and let 
$ \mathcal{E}: 0 \rightarrow K \rightarrow E \rightarrow Q \rightarrow 0 $ be
an extension of separated $ G $-algebras with bounded linear splitting. Then there are two natural exact sequences
\begin{equation*}
\xymatrix{
HL^G_0(A,K)\; \ar@{->}[r] \ar@{<-}[d] & HL^G_0(A,E) \ar@{->}[r] & HL^G_0(A,Q) \ar@{->}[d] \\
HL^G_1(A,Q)\; \ar@{<-}[r] & HL^G_1(A,E) \ar@{<-}[r] & HL^G_1(A,K)
}
\end{equation*}
and
\begin{equation*}
\xymatrix{
HL^G_0(Q,A)\; \ar@{->}[r] \ar@{<-}[d] & HL^G_0(E,A) \ar@{->}[r] & HL^G_0(K,A) \ar@{->}[d] \\
HL^G_1(K,A)\; \ar@{<-}[r] & HL^G_1(E,A) \ar@{<-}[r] & HL^G_1(Q,A)
}
\end{equation*}
The horizontal maps in these diagrams are induced by the maps 
in $ \mathcal{E} $ and the vertical maps are, up to a sign, given by composition product 
with an element $ \ch(\mathcal{E}) $ in $ HL^G_1(Q,K) $ naturally associated to the extension. 
Analogous statements hold for the analytic theory. 
\end{theorem}
Upon tensoring the given extension $ \mathcal{E} $ with $ \mathcal{K}_G $ we obtain 
an extension of separated $ G $-algebras with equivariant bounded linear splitting. As
in \cite{Voigtepch} we may suppress the algebra $ \mathcal{K}_G $ from our notation and assume 
that we are given an extension
\begin{equation*}\label{eqext}
\xymatrix{
K\; \ar@{>->}[r]^{\iota} & E \ar@{->>}[r]^{\pi} & Q 
}
\end{equation*}
of separated $ G $-algebras together with an equivariant bounded linear splitting $ \sigma: Q \rightarrow E $ 
for the quotient map $ \pi: E \rightarrow Q $. \\
We denote by $ X_G(\mathcal{T}E:\mathcal{T}Q) $ the kernel of the map $ X_G(\mathcal{T}\pi):
X_G(\mathcal{T}E) \rightarrow X_G(\mathcal{T}Q)) $ induced by $ \pi $. The splitting $ \sigma $ 
yields a direct sum decomposition $ X_G(\mathcal{T}E) = X_G(\mathcal{T}E:\mathcal{T}Q) \oplus X_G(\mathcal{T}Q) $ 
of $ \ayd $-modules. Moreover there is a natural chain map 
$ \rho: X_G(\mathcal{T}K) \rightarrow X_G(\mathcal{T}E:\mathcal{T}Q) $. \\
Theorem \ref{Excision} is a consequence of the following result. 
\begin{theorem}\label{Excision2} The map $ \rho: X_G(\mathcal{T}K) \rightarrow
X_G(\mathcal{T}E:\mathcal{T}Q) $ is a homotopy equivalence.
\end{theorem}
\proof The proof follows the arguments given in \cite{Meyerthesis}, \cite{Voigtepch}. Let $ \mathfrak{L} \subset \mathcal{T}E $ be the left ideal generated by
$ K \subset \mathcal{T}E $. Then $ \mathfrak{L} $ is a separated $ G $-algebra and we obtain 
an extension 
\begin{equation*}
\xymatrix{
   N\; \ar@{>->}[r] & \mathfrak{L} \ar@{->>}[r]^{\tau} & K \\
 }
\end{equation*}
of separated $ G $-algebras where $ \tau: \mathfrak{L} \rightarrow K $ is induced by the canonical 
projection $ \tau_E: \mathcal{T}E \rightarrow E $. 
As in \cite{Voigtepch} one shows that the inclusion 
$ \mathfrak{L} \subset \mathcal{T}E $ induces a homotopy equivalence 
$ \psi: X_G(\mathfrak{L}) \rightarrow X_G(\mathcal{T}E:\mathcal{T}Q) $. 
The inclusion $ \mathcal{T}K \rightarrow \mathfrak{L} $ induces a morphism of extensions
from $  0 \rightarrow \mathcal{J}K \rightarrow \mathcal{T}K \rightarrow K \rightarrow 0 $ 
to $ 0 \rightarrow N \rightarrow \mathfrak{L} \rightarrow K \rightarrow 0 $. 
The algebra $ N $ is analytically nilpotent and the splitting 
homomorphism $ v: \mathfrak{L} \rightarrow \mathcal{T} \mathfrak{L} $ for the canonical projection constructed by Meyer in 
\cite{Meyerthesis} is easily seen to be equivariant. Using homotopy invariance it follows that the induced chain map $ X_G(\mathcal{T}K) \rightarrow 
X_G(\mathfrak{L}) $ is a homotopy equivalence. This yields the assertion. \qed 

\section{Comparison between analytic and local cyclic homology} \label{seccomp}

In this section we study the relation between equivariant analytic and local cyclic homology. We exhibit a 
special case in which the analytic and local theories agree. This allows to do some elementary  
calculations in equivariant local cyclic homology. Our discussion follows closely the treatment 
by Meyer in \cite{Meyerlocalcyclic}. For the convenience of the reader we reproduce some of his results. \\
A bornological vector space $ V $ is called subcomplete 
if the canonical map $ V \rightarrow V^c $ is a bornological embedding with locally dense range. 
\begin{prop} \label{subcompletechar}
Let $ V $ be a separated bornological vector space. The following conditions are equivalent: 
\begin{bnum}
\item[a)] $ V $ is subcomplete. 
\item[b)] for every small disk $ S \subset V $ there is a small disk $ T \subset V $ containing $ S $ such that 
every $ S $-Cauchy sequence that converges in $ V $ is already $ T $-convergent. 
\item[c)] for every small disk $ S \subset V $ there is a small disk $ T \subset V $ containing $ S $ such that 
every $ S $-Cauchy sequence which is a null sequence in $ V $ is already a $ T $-null sequence.
\item[d)] for every small disk $ S \subset V $ there is a small disk $ T \subset V $ containing $ S $ such that 
$$
\ker(\bra S \ket^c \rightarrow \bra T \ket^c) = \ker(\bra S \ket^c \rightarrow \bra U \ket^c)
$$
for all small disks $ U $ containing $ T $. 
\item[e)] for every small disk $ S \subset V $ there is a small disk $ T \subset V $ containing $ S $ such that 
$$
\ker(\bra S \ket^c \rightarrow \bra T \ket^c) = \ker(\bra S \ket^c \rightarrow V^c).
$$ 
\end{bnum}
\end{prop}
\proof $ a) \Rightarrow b) $ Let $ S \subset V $ be a small disk. Then there exists a small disk $ R \subset V^c $ such that 
every $ S $-Cauchy sequence is $ R $-convergent. Since $ V \rightarrow V^c $ is a bornological embedding 
the disk $ T = R \cap V $ is small in $ V $. By construction, every $ S $-Cauchy sequence that converges in 
$ V $ is already $ T $-convergent. 
$ b) \Rightarrow c) $ is clear since $ V $ is separated. 
$ c) \Leftrightarrow d) $ Let $ U $ be a small disk containing $ S $. Then the kernel of the map $ \bra S \ket^c \rightarrow \bra U \ket^c $ 
consists of all $ S $-Cauchy sequences which are $ U $-null sequences. Since a null sequence in $ V $ is a null sequence in $ U $ for 
some small disk $ U $ the claim follows. $ d) \Rightarrow e) $ Let $ \Analysis(V) $ be the inductive system of normed spaces obtained as the 
dissection of the bornological vector space $ V $. Condition $ d) $ implies that the direct limit of $ \Analysis(V)^c $ is 
automatically separated. That is, $ V^c = \varinjlim \Analysis(V)^c $ is equal to the vector space direct limit of the system $ \Analysis(V)^c $ 
with the quotient bornology. Hence $ \ker(\bra S \ket \rightarrow V^c) = \bigcup \ker(\bra S \ket^c \rightarrow \bra U \ket^c) $ where the union is taken over 
all small disks $ U $ containing $ S $. $ e) \Rightarrow a) $ For each small disk $ S $ in $ V $ let us define 
$ \bra \bra S \ket \ket = \bra S \ket^c/\ker(\bra S \ket^c \rightarrow V^c) $. According to $ e) $ the resulting inductive system is isomorphic to $ \Analysis(V)^c $ 
and $ \varinjlim \bra \bra S \ket \ket \cong V^c $. Assume that $ x \in \ker(\bra S \ket \rightarrow \bra \bra S \ket \ket) $. Then 
$ x \in \ker(\bra S \ket \rightarrow \bra T \ket^c) $ for some small disk $ T $ containing 
$ S $. This implies $ x = 0 $ since the maps $ \bra S \ket \rightarrow \bra T \ket \rightarrow \bra T \ket^c $ are 
injective. Hence $ \bra S \ket \rightarrow \bra \bra S \ket \ket $ is injective for all small disks $ S $. It follows that $ \iota: V \rightarrow V^c $ 
is a bornological embedding with locally dense range. \qed \\
We are interested in conditions which imply that the space $ \Omega(A) $ for a separated 
bornological algebra $ A $ is subcomplete. As usual, we consider $ \Omega(A) $ as a bornological 
vector space with the analytic bornology. Given a small set $ S \subset A $ we shall write $ \Omega(S) $ for the disked hull of 
$$
S \cup \bigcup_{n = 1}^\infty S^{\otimes n + 1} \cup S^{\otimes n} 
$$
inside $ \Omega(A) $ where we use the canonical identification $ \Omega^n(A) = A^{\otimes n + 1} \oplus A^{\otimes n} $ 
for the space of differential forms. Remark that the sets $ \Omega(S) $ generate the analytic bornology. 
\begin{definition}
A separated bornological algebra $ A $ is called tensor subcomplete if the space $ \Omega(A) $ is subcomplete. 
\end{definition}
Let us call the tensor powers $ V^{\otimes n} $ for $ n \in \mathbb{N} $ of a bornological vector space $ V $ uniformly 
subcomplete provided the following condition is satisfied. For every small disk $ S \subset V $ there is a 
small disk $ T \subset V $ containing $ S $ such that, independent of $ n \in \mathbb{N} $, any $ S^{\otimes n} $-Cauchy sequence 
which is a null sequence in $ V^{\otimes n} $ is already a $ T^{\otimes n} $-null sequence. In particular, the spaces $ V^{\otimes n} $ are 
subcomplete for all $ n $ in this case. 
\begin{lemma}\label{subcomplemma}
A separated bornological algebra $ A $ is tensor subcomplete iff the tensor powers $ A^{\otimes n} $ for 
$ n \in \mathbb{N} $ are uniformly subcomplete. 
\end{lemma}
\proof Assume first that the space $ \Omega(A) $ is subcomplete. Let $ S \subset A $ be a small disk and let $ (x_k)_{k \in \mathbb{N}} $ be a 
$ S^{\otimes n} $-Cauchy sequence which is a null sequence in $ A^{\otimes n} $. 
We write $ i_n: A^{\otimes n} \rightarrow \Omega(A) $ and $ p_n: \Omega(A) \rightarrow A^{\otimes n} $ 
for the natural inclusion and projection onto one of the direct summands $ A^{\otimes n} $ in $ \Omega(A) $. 
The maps $ i_n $ and $ p_n $ are clearly bounded. In particular, the image of $ (x_k)_{k \in \mathbb{N}} $ 
under $ i_n $ is a $ \Omega(S) $-Cauchy sequence which is a null sequence in $ \Omega(A) $. Hence it is a $ \Omega(T) $-null sequence 
for some $ T \subset A $. Since $ p_n(\Omega(T)) = T^{\otimes n} $ and $ x_k = p_n i_n(x_k) $ it follows that 
the sequence $ (x_k)_{k \in \mathbb{N}} $ is a $ T^{\otimes n} $-null sequence. Moreover the choice of $ T $ does not depend on $ n $.  
This shows that the tensor powers $ A^{\otimes n} $ are uniformly subcomplete. \\
Conversely, assume that the tensor powers $ A^{\otimes n} $ are uniformly subcomplete. Let $ S \subset A $ 
be a small disk and let $ T \subset A $ be a small disk such that $ S^{\otimes n} $-Cauchy which are 
null sequences in $ A^{\otimes n} $ are $ T^{\otimes n} $-null sequences. In addition we may assume 
$ 2 S \subset T $. Let us write $ P_n: \Omega(A) \rightarrow \Omega(A) $ for the natural projection 
onto the direct summand $ \bigoplus_{j = 1}^n A^{\otimes j} \oplus A^{\otimes j} $. Then $ P_n(\Omega(S)) $ 
is contained in $ \Omega(S) $ and the projections $ P_n $ are equibounded. Moreover $ P_n $ converges 
to the identity uniformly on $ \Omega(S) $ since $ \id - P_n $ has norm $ \leq 2^{-n} $ as 
a map from $ \bra \Omega(S) \ket $ into $ \bra \Omega(2S) \ket \subset \bra \Omega(T) \ket $. Now 
let $ (x_k)_{k \in \mathbb{N}} $ be a null sequence in $ \Omega(A) $ which is $ \Omega(S) $-Cauchy. 
The components of $ P_n(x_k) $ are $ S^{\otimes k} $-Cauchy sequences which are null sequences in $ A^{\otimes k} $ and hence 
$ T^{\otimes k} $-null sequences by hypothesis. Hence $ (P_n(x_k))_{n \in \mathbb{N}} $ is a $ \Omega(T) $-null 
sequence for all $ n $. Moreover $ (P_n(x_k))_{n \in \mathbb{N}} $ converges to $ x_k $ 
for every $ k $, and this convergence is uniform in $ k $. If follows that 
$ (x_k)_{n \in \mathbb{N}} $ is a $ \Omega(T)$-null sequence, and we deduce that $ \Omega(A) $ is tensor subcomplete. \qed \\
Our next aim is to exhibit certain analytical conditions which are sufficient for tensor subcompleteness. 
Recall that a subset $ S $ of a complete bornological vector space $ V $ is called (relatively) compact if 
it is a (relatively) compact subset of the Banach space $ \bra T \ket $ for some small completant disk $ T \subset V $. 
A complete bornological vector space $ V $ is a Schwartz space if every small subset of $ V $ is relatively compact. Every 
Fr\'echet space with the precompact bornology is a Schwartz space. \\
Let $ V $ be a normed space and let $ W $ be an arbitrary bornological vector space. By definition, a sequence $ (f_n)_{n \in \mathbb{N}} $ of bounded linear maps 
$ f_n: V \rightarrow W $ converges uniformly to $ f: V \rightarrow W $ if there exists a small disk $ T \subset W $ 
such that all $ f_n $ and $ f $ are bounded linear maps $ V \rightarrow \bra T \ket $ and the sequence $ (f_n)_{n \in \mathbb{N}} $ 
converges to $ f $ in $ \Hom(V, \bra T \ket) $ in operator norm. 
A bounded linear map $ f: V \rightarrow W $ can be approximated uniformly on compact subsets by 
finite rank operators if for every compact disk $ S \subset V $ there exists a sequence $ (f_n)_{n \in \mathbb{N}} $ of 
finite rank operators $ f_n: V \rightarrow W $ such that $ f_n $ converges uniformly to $ f $ in $ \Hom(\bra S \ket, W) $.  
An operator $ f: V \rightarrow W $ is of finite rank if it is contained in the image of 
the natural map from the uncompleted tensor product $ W \otimes V' $ into $ \Hom(V,W) $ 
where $ V' = \Hom(V, \mathbb{C}) $ is the dual space of $ V $. 
By definition, a complete bornological vector space $ V $ satisfies the (global) approximation property if the identity map 
on $ V $ can be approximated uniformly on compact subsets by finite rank operators. \\
We recall that a bornological vector space $ V $ is regular if the bounded linear functionals on $ V $ separate the points of $ V $. 
Let us remark that there is also a local version of the approximation property which is equivalent to the global one if we restrict 
attention to regular spaces. 
Finally, we point out that for a Fr\'echet space with the precompact bornology the bornological approximation property is equivalent to Grothendieck's approximation 
property \cite{Meyerborntop}. 
\begin{prop}\label{apptensorsub}
Let $ A $ be a bornological algebra whose underlying bornological vector space is a Schwartz space satisfying the approximation property. Then $ A $ is tensor subcomplete. 
\end{prop}
\proof According to lemma \ref{subcomplemma} it suffices to show that the tensor powers of $ A $ are uniformly subcomplete. 
Let $ S \subset A $ be a small disk. We may assume $ S $ is compact and that there is a completant small disk $ T \subset A $ 
containing $ S $ such that the inclusion $ \bra S \ket \rightarrow \bra T \ket $ can be approximated 
uniformly by finite rank operators on $ A $. We will show that $ \ker(\bra S \ket^{\cotimes n} \rightarrow \bra U \ket^{\cotimes n}) =
\ker(\bra S \ket^{\cotimes n} \rightarrow \bra T \ket^{\cotimes n}) $ 
for every completant small disk $ U $ containing $ T $. As in the proof of proposition \ref{subcompletechar} this statement easily 
implies that the tensor powers $ A^{\otimes n} $ are uniformly subcomplete. \\
Take an element $ x \in \ker(\bra S \ket^{\cotimes n} \rightarrow \bra U \ket^{\cotimes n}) $. Then there is a compact 
disk $ K \subset \bra S \ket $ such that $ x \in K^{\cotimes n} $. Since $ A $ is regular we find a sequence $ (f_k)_{k \in \mathbb{N}} $ 
of finite rank operators $ f_k: A \rightarrow \bra T \ket $ approximating the inclusion map uniformly on $ K $. 
The uniform convergence of the operators $ f_k $ on $ K $ implies that $ f^{\cotimes n}_k $ converges uniformly towards the 
canonical map $ \bra K \ket^{\cotimes n} \rightarrow \bra T \ket^{\cotimes n} $. In particular, the image 
of $ x $ in $ \bra T \ket^{\cotimes n} $ is the limit of $ f^{\cotimes n}_k(x) $. 
Since the finite rank maps $ f_k $ are restrictions of maps defined on $ \bra U \ket $ and $ x $ is in 
the kernel of $ \bra S \ket^{\cotimes n} \rightarrow \bra U \ket^{\cotimes n} $ we have 
$ f^{\cotimes n}_k(x) = 0 $ for all $ k $. Hence $ x \in \ker(\bra S \ket^{\cotimes n} \rightarrow \bra T \ket^{\cotimes n}) $ 
as desired. \qed \\
It follows in particular that the algebra $ A \otimes \mathcal{K}_G $ is tensor subcomplete 
provided $ A $ is a Schwartz space satisfying the approximation property. 
\begin{prop}\label{limitacyclic}
Let $ G $ be a totally disconnected group and let $ A $ be a $ G $-algebra whose underlying bornological vector space is a Schwartz space 
satisfying the approximation property. Then the canonical chain map 
$$ 
X_G(\mathcal{T}(A \otimes \mathcal{K}_G))^{\mathbb{L}c} \rightarrow X_G(\mathcal{T}(A \otimes \mathcal{K}_G))^c 
$$ 
induces an isomorphism in the local derived category. 
\end{prop}
\proof Let us abbreviate $ X = X_G(\mathcal{T}(A \otimes \mathcal{K}_G)) $ and remark that the $ \ayd $-module $ X $ can be written in the form 
$ X = \AYD(G) \otimes V $ for a separated bornological vector space $ V $. Using this observation and proposition \ref{apptensorsub} one checks easily that 
the inductive system $ \Analysis(X)^c $ is essentially injective. Due to proposition \ref{holimlocproj} 
it follows that the natural map $ X^{\mathbb{L}c} \cong \holim (\Analysis(X)^c) \rightarrow \varinjlim (\Analysis(X)^c) = X^c $ 
is a local homotopy equivalence. This yields the claim. \qed \\
An analogous argument shows that $ X_G(\mathcal{T}\mathbb{C})^{\mathbb{L}c} \rightarrow X_G(\mathcal{T}\mathbb{C})^c  $ 
is a local homotopy equivalence. It follows that there is a chain of canonical isomorphisms
$$ 
X_G(\mathcal{T}\mathbb{C})^{\mathbb{L}c} \cong X_G(\mathcal{T}\mathbb{C})^c\cong X_G(\mathbb{C}) = \mathcal{O}_G[0] 
$$ 
in the local derived category. In fact, the seond isomorphism is a consequence of the fact that $ \mathbb{C} $ is analytically quasifree 
combined with homotopy invariance. The last equality is established in \cite{Voigtepch}. \\
Consider in particular the case that $ G $ is a compact group. Then the paracomplex $ \mathcal{O}_G[0] $ is primitive. 
Taking into account stability, this yields
$$
HL^G_*(\mathbb{C}, B) = H_*(\SHom_G(\mathcal{O}_G[0], X_G(\mathcal{T}B)^{\mathbb{L}c})), 
$$
and analogously we have 
$$
HA^G_*(\mathbb{C}, B) = H_*(\SHom_G(\mathcal{O}_G[0], X_G(\mathcal{T}B)^c))
$$
for every $ G $-algebra $ B $. We conclude that there exists a natural transformation $ HL^G_*(\mathbb{C},B) \rightarrow HA^G_*(\mathbb{C},B) $ 
between equivariant local and analytic cyclic homology if the group is compact.
\begin{prop}\label{hlhacompcomp}
Let $ G $ be compact and let $ B $ be a $ G $-algebra whose underlying bornological vector space is a Schwartz space satisfying the approximation property. 
Then the natural map 
$$
HL^G_*(\mathbb{C}, B) \rightarrow HA^G_*(\mathbb{C}, B)
$$
is an isomorphism. In particular, there is a canonical isomorphism 
$$
HL^G_*(\mathbb{C}, \mathbb{C}) \cong HA^G_*(\mathbb{C}, \mathbb{C}) = \mathcal{R}(G)
$$
where $ \mathcal{R}(G) $ is the algebra of conjugation invariant smooth functions on $ G $. 
\end{prop}
\proof Using stability, the first assertion follows from proposition \ref{limitacyclic} and the fact that $ \mathcal{O}_G[0] $ is primitive. 
For the second claim observe that $ \mathcal{R}(G) = (\mathcal{O}_G)^G $ is the invariant part 
of $ \mathcal{O}_G $. \qed 

\section{The isoradial subalgebra theorem}\label{seciso}

In this section we discuss the notion of an isoradial subalgebra and prove the isoradial subalgebra theorem 
which states that equivariant local cyclic homology is invariant under the passage to 
isoradial subalgebras. \\
Recall that a subset $ \mathcal{V} $ of a bornological vector space $ V $ is called 
locally dense if for any small subset $ S \subset V $ there exists a small disk $ T \subset V $ 
such that any $ v \in S $ is the limit of a $ T $-convergent 
sequence with entries in $ \mathcal{V} \cap T $. Moreover recall that a separated 
(complete) bornological algebra $ A $ is locally multiplicative iff 
it is isomorphic to an inductive limit of (complete) normed algebras. The following definition 
is taken from \cite{Meyerborntop}. 
\begin{definition}\label{isosubalg}
Let $ \A $ and $ A $ be complete locally multiplicative bornological algebra. A bounded homomorphism $ \iota: \mathcal{A} \rightarrow A $ between bornological algebras is called isoradial if it has locally dense range and 
\begin{equation*}
\rho(\iota(S); A) = \rho(S;\mathcal{A})
\end{equation*}
for all small subsets $ S \subset \mathcal{A} $. If in addition $ \iota $ is injective we say that 
$ \mathcal{A} $ is an isoradial subalgebra of $ A $.
\end{definition}
We will frequently identify $ \mathcal{A} $ with its image $ \iota(\mathcal{A}) \subset A $ 
provided $ \iota: \mathcal{A} \rightarrow A $ is an injective bounded homomorphism.
However, note that the bornology of $ \mathcal{A} $ is usually finer than the subspace bornology 
on $ \iota(\mathcal{A}) $. Remark in addition that the inequality
$ \rho(\iota(S);A) \leq \rho(S; \mathcal{A})$ is automatic for every small subset $ S \subset \mathcal{A} $. 
If $ \A $ and $ A $ are $ G $-algebras and $ \iota: \A \rightarrow A $ is an 
equivariant homomorphism defining an isoradial subalgebra we say that $ \mathcal{A} $ is an isoradial $ G $-subalgebra of $ A $ . \\
Assume that $ \iota: \mathcal{A} \rightarrow A $ is an equivariant homomorphism and consider the equivariant homomorphism
$ i: \mathcal{A} \otimes \mathcal{K}_G \rightarrow A \otimes \mathcal{K}_G $ obtained by tensoring 
$ \iota $ with the identity map on $ \mathcal{K}_G $. 
It is shown in \cite{Meyerborntop} that isoradial homomorphisms are preserved under tensoring with nuclear locally 
multiplicative algebras. In particular, this yields the following statement. 
\begin{prop}\label{isotensor} If $ \iota: \mathcal{A} \rightarrow A $ is an isoradial $ G $-subalgebra then 
$ i: \mathcal{A} \otimes \mathcal{K}_G \rightarrow A \otimes \mathcal{K}_G $ is an 
isoradial $ G $-subalgebra as well. 
\end{prop}
Note that the algebra $ \mathcal{K}_G $ carries the fine bornology which implies that tensor products of $ \mathcal{K}_G $ with 
complete spaces are automatically complete. \\
Let us now formulate and prove the isoradial subalgebra theorem. 
\begin{theorem} \label{isosubalgebratheorem1} 
Let $ \iota: \mathcal{A} \rightarrow A $ be an isoradial 
$ G $-subalgebra. Suppose that there exists 
a sequence $ (\sigma_n)_{n \in I} $ of bounded linear maps 
$ \sigma_n: A \rightarrow \mathcal{A} $ such that 
for each completant small disk $ S \subset $ A the maps 
$ \iota \sigma_n $ converge uniformly towards the 
inclusion map $ \bra S \ket \rightarrow A $. Then the class 
$ [\iota] \in HL^G_*(\mathcal{A},A) $ is invertible. 
\end{theorem}
Note that the existence of bounded linear maps $ \sigma_n: A \rightarrow \A $ with these properties already implies that 
$ \A \subset A $ is locally dense. We point out that the maps $ \sigma_n $ in theorem \ref{isosubalgebratheorem1} are not assumed to 
be equivariant. \\
In fact, as a first step in the proof we shall modify these maps in order to 
obtain equivariant approximations. Explicitly, let us define equivariant bounded linear maps 
$ s_n: A \otimes \mathcal{K}_G \rightarrow \mathcal{A} \otimes \mathcal{K}_G $ by 
\begin{equation*}
s_n(a \otimes k)(r, t) = 
t \cdot \sigma_n(t^{-1} \cdot a) k(r,t)
\end{equation*}
where we view elements in $ \mathcal{A} \otimes \mathcal{K}_G $ as smooth function on $ G \times G $ 
with values in $ \mathcal{A} $. As above we write $ i $ for the equivariant homomorphism $ \mathcal{A} \otimes \mathcal{K}_G 
\rightarrow A \otimes \mathcal{K}_G $ induced by $ \iota $. 
Since the maps $ \iota \sigma_n $ converge to the identity uniformly on small subsets of $ A $ by assumption, 
the maps $ i s_n $ converge to the identity uniformly on small subsets of $ A \otimes \mathcal{K}_G $. \\
We deduce that theorem \ref{isosubalgebratheorem1} is a consequence of the following theorem. 
\begin{theorem} \label{isosubalgebratheorem2} Let $ \iota: \mathcal{A} \rightarrow A $ be an isoradial 
$ G $-subalgebra. Suppose that there exists 
a sequence $ (\sigma_n)_{n \in I} $ of equivariant bounded linear maps 
$ \sigma_n: A \rightarrow \mathcal{A} $ such that 
for each completant small disk $ S \subset $ A the maps 
$ \iota \sigma_n $ converge uniformly towards the 
inclusion map $ \bra S \ket \rightarrow A $. Then the chain map $ X_G(\mathcal{T}\mathcal{A}) \rightarrow X_G(\mathcal{T}A) $ 
induced by $ \iota $ is a local homotopy equivalence. 
\end{theorem}
The proof of theorem \ref{isosubalgebratheorem2} is divided into several steps. 
Let $ S \subset A $ be a small completant multiplicatively closed disk. 
By the definition of uniform convergence, there exists a small completant disk $ T \subset A $ 
containing $ S $ such that $ \iota \sigma_n $ defines a bounded linear map $ \bra S \ket \rightarrow \bra T \ket $ for every $ n $ and 
the sequence $ (\iota \sigma_n)_{n \in \mathbb{N}} $ converges to the natural inclusion map in $ \Hom(\bra S \ket, \bra T \ket) $ in 
operator norm. Hence there exists a null sequence
$ (\epsilon_n)_{n \in \mathbb{N}} $ of positive real numbers such that $ \iota \sigma_n(x) - x \in \epsilon_n T $ for 
all $ x \in S $. After rescaling with a positive scalar $ \lambda $ we may assume that $ T $ is multiplicatively 
closed and that $ S \subset \lambda T $. Using the formula 
\begin{align*}
\omega&_{\iota \sigma_n}(x,y) = \iota \sigma_n(xy) - \iota \sigma_n(x)\iota \sigma_n(y) \\
&= (\iota \sigma_n(xy) - xy) - (\iota \sigma_n(x) - x)(\iota \sigma_n(y) - y) - x (\iota \sigma_n(y) - y) - (\iota \sigma_n(x) - x) y 
\end{align*}
for $ x,y \in S $ and that $ T $ is multiplicatively closed we see that for any given 
$ \epsilon > 0 $ we find $ N \in \mathbb{N} $ such that  $ \omega_{\iota \sigma_n}(S,S) \subset \epsilon T $ for $ n \geq N $. 
Remark that we have $ \omega_{\iota \sigma_n} = \iota \omega_{\sigma_n} $ since $ \iota $ is a homomorphism. We deduce 
$$
\lim_{n \rightarrow \infty} \rho(\iota \omega_{\sigma_n}(S,S); A) = 0
$$
using again that $ T $ is multiplicatively closed. This in turn implies 
$$
\lim_{n \rightarrow \infty} \rho(\omega_{\sigma_n}(S,S); \A) = 0
$$
since $ \A \subset A $ is an isoradial subalgebra. 
This estimate will be used to obtain local inverses for the chain map induced by $ \iota $. \\
We need some preparations. Let $ B $ and $ C $ be arbitrary separated $ G $-algebras. Any equivariant bounded linear map $ f: B \rightarrow C $ extends to an equivariant homomorphism 
$ \mathcal{T}f: \mathcal{T}B \rightarrow \mathcal{T} C $. This homomorphism is bounded 
iff $ f $ has analytically nilpotent curvature. 
\begin{lemma} \label{taupreservespecr}
Let $ C $ be a separated bornological algebra and let $ S \subset \mathcal{T}C $ be small. Then 
$$
\rho(\tau_C(S); C) = \rho(S; \mathcal{T}C)
$$
where $ \tau_C: \mathcal{T}C \rightarrow C $ is the quotient homomorphism. 
\end{lemma}
\proof Taking into account that $ \tau_C $ is a bounded homomorphism it suffices to show that $ \rho(\tau_C(S); C) < 1 $ implies $ \rho(S; \mathcal{T}C) \leq 1 $. 
We may assume that the set $ S $ is of the form 
$$
S = R + [T](dT dT)^\infty
$$
where $ R \subset C $ and $ T \subset C $ are small disks. If $ \rho(\tau_C(S); C)< 1 $ we find $ \lambda > 1 $ such that 
$ (\lambda R)^\infty \subset C $ is small. Let us choose $ \mu $ such that $ \lambda^{-1} + \mu^{-1} < 1 $ and 
consider the small disk $ P = \mu (\lambda R)^\infty $ in $ C $. By construction we have 
$ R \cdot [P] \subset \lambda^{-1} P $ as well as $ dR d[P] \subset \mu^{-1} dP dP $. 
Moreover the disked hull $ I $ of 
$ P \cup [P] (dP dP)^\infty $ 
is a small subset of $ \mathcal{T}C $ which contains $ R $. Now consider $ x \in R $ and 
$ [y_0] dy_1 \cdots dy_{2n} \in [P](dPdP)^n $. Since
$$
x \circ [y_0] dy_1 \cdots dy_{2n} = x [y_0] dy_1 \cdots dy_{2n} + dx d[y_0] dy_1 \cdots dy_{2n}
$$
the previous relations yield $ \nu R \circ I \subset I $ for some $ \nu > 1 $. By induction we see that the multiplicative closure $ Q $ 
of $ \nu R $ in $ \mathcal{T}C $ is small. 
Choose $ \eta $ such that $ \nu^{-1} + 2 \eta^{-1} < 1 $, set 
$$ 
K = \eta [T] (dT dT)^\infty 
$$ 
and let $ L $ be the multiplicative closure of $ [Q] \circ K \circ [Q] $. 
By construction, the set $ [Q] \circ K \circ [Q] $ is contained in the analytically nilpotent algebra $ \mathcal{J}C $ 
which implies that $ L \subset \mathcal{T}C $ is small. Let $ J \subset \mathcal{T}C $ be the disked hull of the set $ Q + L $. Then 
$ J $ is small and we have $ S \subset J $. In addition, it is straightforward to check $ R \circ J \subset \nu^{-1} J $ and 
$ [T](dT dT)^\infty \circ J \subset 2 \eta^{-1} J $ which shows $ S \circ J \subset J $. In the same way as above 
it follows that $ S^\infty \subset \mathcal{T}C $ is small and deduce $ \rho(S; \mathcal{T}C) \leq 1 $. \qed 
\begin{lemma}\label{localbounded}
Let $ f: B \rightarrow C $ be an equivariant bounded linear map between separated $ G $-algebras. Consider the 
induced chain map $ X_G(\mathcal{T}f): X_G(\mathcal{T}B) \rightarrow X_G(\mathcal{T}C) $. Given a small subset $ S \subset 
X_G(\mathcal{T}B) $ there exists a small subset $ T \subset B $ such that 
$ X_G(\mathcal{T}f) $ is bounded on the primitive submodule generated by $ S $ provided 
$ \omega_f(T,T)^\infty $ is small. 
\end{lemma} 
\proof It suffices to show that, given a small set $ S \subset \mathcal{T}B $, there exists 
a small set $ T \subset B $ such that $ \mathcal{T}f(S) \subset\mathcal{T}C $ is 
small provided $ \omega_f(T,T)^\infty $ is small. 
We may assume that $ S $ is of the form $ [R] (dR dR)^\infty $ for some small set $ R \subset B $. 
Let $ F: B \rightarrow \mathcal{T}C $ be the bounded linear map obtained by composing $ f $ with the 
canonical bounded linear splitting $ \sigma_C: C \rightarrow \mathcal{T}C $. The homomorphism $ \mathcal{T}f: \mathcal{T}B \rightarrow \mathcal{T}C $ 
is given by 
$$
\mathcal{T}f([x_0] dx_1 \cdots dx_{2n}) = [F(x_0)] \omega_F(x_1,x_2) \cdots \omega_F(x_{2n - 1}, x_{2n})
$$
which shows that $ \mathcal{T}f(S) $ is small provided $ \omega_F(R,R)^\infty $ is small. 
Consider the natural projection $ \tau_C: \mathcal{T}C \rightarrow C $. According to lemma \ref{taupreservespecr} the homomorphism $ \tau_C $ 
preserves the spectral radii of all small subsets in $ \mathcal{T}C $. Using $ \tau_C(\omega_F(R,R)) = \omega_f(R,R) $ we see that 
$ \omega_F(R,R)^\infty $ is small provided $ \rho(\omega_f(R,R)) < 1 $. 
Setting $ T = \lambda R $ for some $ \lambda > 1 $ yields the assertion. \qed \\
Let us come back to the proof of theorem \ref{isosubalgebratheorem2}. If $ P \subset X_G(\mathcal{T}A) $ is a primitive subparacomplex 
then lemma \ref{localbounded} shows that $ \sigma_n $ induces a bounded chain map $ P \rightarrow X_G(\mathcal{T}\A) $ provided 
$ n $ is large enough. In fact, we will prove that the maps $ \sigma_n $ can be used to define bounded local homotopy inverses to the chain map 
$ \iota_*: X_G(\mathcal{T}\A) \rightarrow X_G(\mathcal{T}A)$ induced by $ \iota $. \\
More precisely, let $ k_n: A \rightarrow A \otimes \mathbb{C}[t] $ be the equivariant 
bounded linear map given by $ k_n(x)(t) = (1 - t)(\iota \sigma_n)(x) + t x $. Here $ \mathbb{C}[t] $ is equipped with 
the bornology induced from $ C^\infty[0,1] $. Observe that the maps $ k_n $ converge to the homomorphism sending $ x $ to $ x \otimes 1 $ 
uniformly on small subsets of $ A $. The same reasoning as for the maps $ \iota \sigma_n $ above shows
$$
\lim_{n \rightarrow \infty} \rho(\omega_{k_n}(S,S); A \otimes \mathbb{C}[t]) = 0
$$
for all small subsets $ S \subset A $. Now assume that $ P \subset X_G(\mathcal{T}A) $ 
is a primitive subparacomplex. 
According to lemma \ref{localbounded} there exists $ N \in \mathbb{N} $ such that the induced chain map 
$ X_G(\mathcal{T}A) \rightarrow X_G(\mathcal{T}(A \otimes \mathbb{C}[t])) $ 
is bounded on $ P $ for all $ n > N $. We compose this map with the 
chain homotopy between the evaluation maps at $ 0 $ and $ 1 $ arising from homotopy invariance 
to get a bounded 
$ \ayd $-map $ K_n: P \rightarrow X_G(\mathcal{T}A) $ of degree 
one which satisfies $ \partial K_n + K_n \partial = \id - (\iota \sigma_n)_* $ on $ P $. \\
Similarly, consider the equivariant bounded linear map 
$ h_n: \A \rightarrow \A \otimes \mathbb{C}[t] $ given by $ h_n(x)(t) = (1 - t)(\sigma_n \iota)(x) + t x $ 
and observe that $ (\iota \otimes \id) h_n = k_n \iota $. Since the algebra $ C^\infty[0,1] $ is nuclear 
the inclusion $ \A \otimes \mathbb{C}[t] \rightarrow A \otimes \mathbb{C}[t] $ preserves the spectral radii of 
small subsets \cite{Meyerborntop}. Hence the above spectral radius estimate for $ k_n $ implies 
$$
\lim_{n \rightarrow \infty} \rho(\omega_{h_n}(S,S); \A \otimes \mathbb{C}[t]) = 0
$$
for all small subsets $ S \subset \A $. Now let $ Q \subset X_G(\mathcal{T}\A) $ be 
a primitive subparacomplex. For $ n $ sufficiently large we obtain in the same way as above a bounded $ \ayd $-map 
$ H_n: Q \rightarrow X_G(\mathcal{T}\A) $ of degree $ 1 $ such that 
$ \partial H_n + H_n \partial = \id - (\sigma_n \iota)_* $ on $ Q $. \\
Using these considerations it is easy to construct bounded local contracting homotopies for the 
mapping cone of the chain map $ \iota_*: X_G(\mathcal{T}\A) \rightarrow X_G(\mathcal{T}A) $. This shows that 
$ \iota_* $ is a local homotopy equivalence and completes the proof of theorem \ref{isosubalgebratheorem2}. 

\section{Applications of the isoradial subalgebra theorem} \label{secsmoothing}

In this section we study some consequences of the isoradial subalgebra theorem in connection with 
$ C^* $-algebras. This is needed to show that equivariant local cyclic homology is a continuously and 
$ C^* $-stable functor on the category of $ G $-$ C^* $-algebras. 
Moreover, we discuss isoradial subalgebras arising from regular smooth functions on simplicial complexes \cite{Voigtbs}. \\
In the sequel we write $ A \otimes B $ for the (minimal) tensor product of two 
$ C^* $-algebras $ A $ and $ B $. We will only consider such tensor products when one of 
the involved $ C^* $-algebras is nuclear, hence the $ C^* $-tensor product is in fact uniquely defined in these situations. 
Moreover, our notation should not lead to confusion with the algebraic tensor product since we will not have to 
work with algebraic tensor products of $ C^* $-algebras at all. All $ C^* $-algebras are equipped with the precompact bornology when
they are considered as bornological algebras. \\
As a technical preparation we have to examine how the smoothing of $ G $-$ C^* $-algebras is compatible with 
isoradial homomorphisms. Let us recall from \cite{Meyersmoothrep} that a representation $ \pi $ of $ G $ on a complete 
bornological vector space $ V $ is continuous if the adjoint of $ \pi $ defines a 
bounded linear map $ [\pi]: V \rightarrow C(G,V) $ where $ C(G,V) $ is the space of continuous functions 
on $ G $ with values in $ V $ in the bornological sense. For our purposes it suffices to remark that the representation of 
$ G $ on a $ G $-$ C^* $-algebra equipped with the precompact bornology 
is continuous in the bornological sense. We need the following special cases of results obtained by Meyer in \cite{Meyerborntop}. 
\begin{lemma}\label{isosmooth}
Let $ \mathcal{A} $ and $ A $ be complete locally multiplicative bornological algebras on which $ G $ acts 
continuously. If $ \iota: \mathcal{A} \rightarrow A $ is an equivariant isoradial homomorphism then 
$$ 
\Smooth(\iota): \Smooth(\A) \rightarrow \Smooth(A) 
$$ 
is an isoradial homomorphism as well.
Moreover, if $ C $ is a complete nuclear locally multiplicative $ G $-algebra then the natural homomorphism 
$$ 
\Smooth(A) \cotimes C \rightarrow \Smooth(A \cotimes C) 
$$ 
is isoradial.  
\end{lemma}
\proof It is shown in \cite{Meyerborntop} that the inclusion $ \Smooth(B) \rightarrow  B $ is an isoradial 
subalgebra for every complete locally multiplicative bornological algebra $ B $ on which $ G $ acts continuously. 
This yields easily the first claim. 
In addition, the homomorphism $ \Smooth(A) \cotimes C \rightarrow A \cotimes C $ is isoradial because $ C $ is nuclear \cite{Meyerborntop}. 
Since the action on $ C $ is already smooth it follows that 
$$ 
\Smooth(A) \cotimes C \cong \Smooth(\Smooth(A) \cotimes C) \rightarrow 
\Smooth(A \cotimes C) 
$$ 
is isoradial according to the first part of the lemma. \qed \\ 
Let $ A $ be a  $ G $-$ C^* $-algebra and consider the natural equivariant homomorphism 
$ A \cotimes C^\infty[0,1] \rightarrow C([0,1],A) = A \otimes C[0,1] $. This map induces a bounded equivariant homomorphism
\begin{equation*}
\Smooth(A) \cotimes C^\infty[0,1] \cong \Smooth(A \cotimes C^\infty[0,1]) \rightarrow 
\Smooth(C([0,1], A)) 
\end{equation*}
and we have the following result.  
\begin{prop}\label{homc*} The map $ \Smooth(A) \hat{\otimes} C^\infty[0,1] \rightarrow \Smooth(A \otimes C[0,1]) $ 
is an isoradial $ G $-subalgebra and defines an invertible element in
$$ 
HL^G_0(\Smooth(A) \cotimes C^\infty[0,1], \Smooth(A \otimes C[0,1]))
$$
for every $ G $-$ C^* $-algebra $ A $. 
\end{prop}
\proof It is shown in \cite{Meyerborntop} that the natural inclusion $ \iota: A \cotimes C^\infty[0,1] \rightarrow C([0,1], A) $ is isoradial. 
Hence the homomorphism $ \Smooth(A) \hat{\otimes} C^\infty[0,1] \rightarrow \Smooth(C([0,1],A)) $ 
is isoradial according to lemma \ref{isosmooth}. \\
We choose a family $ \sigma_n: C([0,1], A) \rightarrow A \hat{\otimes} C^\infty[0,1] $ of equivariant smoothing operators 
such that the maps $ \iota \sigma_n $ are uniformly bounded and converge to the identity pointwise. It follows that the 
maps $ \iota \sigma_n $ converge towards the identity uniformly on precompact subsets of $ C([0,1], A) $. 
The maps $ \sigma_n $ induce equivariant bounded linear maps 
$ \sigma_n: \Smooth(C([0,1],A)) \rightarrow \Smooth(A) \cotimes C^\infty[0,1] $ 
satisfying the condition of the isoradial subalgebra theorem \ref{isosubalgebratheorem1}. This yields the assertion. \qed \\
Let $ \mathbb{K}_G = \mathbb{K}(L^2(G)) $ be the algebra of compact operators 
on the Hilbert space $ L^2(G) $. The $ C^* $-algebra $ \mathbb{K}_G $ is equipped with 
the action of $ G $ induced by the regular representation. For every $ G $-$ C^* $-algebra $ A $ we have a natural bounded 
equivariant homomorphism $ A \cotimes \mathcal{K}_G \rightarrow A \otimes \mathbb{K}_G $. 
This gives rise to equivariant homomorphisms
$$
\Smooth(A) \cotimes \mathcal{K}_G \rightarrow \Smooth(A \cotimes \mathcal{K}_G) \rightarrow 
\Smooth(A \otimes \mathbb{K}_G).
$$
Similarly, let $ \mathbb{K} = \mathbb{K}(l^2(\mathbb{N}))$ be the algebra of compact operators on an infinite dimensional 
separable Hilbert space with the trivial $ G $-action. If $ M_\infty(\mathbb{C}) $ denotes the direct limit of the finite 
dimensional matrix algebras $ M_n(\mathbb{C}) $ we have a canonical bounded homomorphism 
$ \Smooth(A) \cotimes M_\infty(\mathbb{C}) \rightarrow \Smooth(A \otimes \mathbb{K}) $. 
\begin{prop}\label{stabkg} 
The homomorphism $ \Smooth(A) \cotimes \mathcal{K}_G \rightarrow 
\Smooth(A \otimes \mathbb{K}_G) $ is an isoradial $ G $-subalgebra and defines an 
invertible element in 
$$ 
HL^G_0(\Smooth(A) \cotimes \mathcal{K}_G, \Smooth(A \otimes \mathbb{K}_G))
$$
for every $ G $-$ C^* $-algebra $ A $. 
An analogous assertion holds for the homomorphism $ \Smooth(A) \cotimes M_\infty(\mathbb{C}) \rightarrow \Smooth(A \otimes \mathbb{K}) $. 
\end{prop}
\proof We will only treat the map $ \Smooth(A) \cotimes \mathcal{K}_G \rightarrow 
\Smooth(A \otimes \mathbb{K}_G) $ since the claim concerning the compact operators with the trivial action is obtained in a similar way. \\
Observe that a small subset of $ \mathcal{K}_G $ is contained in a finite dimensional subalgebra of the form $ M_n(\mathbb{C}) $. Since 
$ A \otimes M_n(\mathbb{C}) $ is a bornological subalgebra of $ A \otimes \mathbb{K}_G $ it follows that the homomorphism 
$ \iota: A \cotimes \mathcal{K}_G \rightarrow A \otimes \mathbb{K}_G $ is isoradial. Due to lemma \ref{isosmooth} the same is 
true for the induced map $ \Smooth(A) \cotimes \mathcal{K}_G \rightarrow \Smooth(A \otimes \mathbb{K}_G) $. 
Since $ G $ is second countable and $ \D(G) \subset L^2(G) $ is dense 
we find a countable orthonormal basis $ (e_n)_{n \in \mathbb{N}} $ of $ L^2(G) $ contained in $ \D(G) $. 
Projecting to the linear subspace $ \mathbb{C}^n \subset L^2(G) $ generated by the vectors $ e_1, \dots, e_n $ defines a 
bounded linear map $ \sigma_n: A \otimes \mathbb{K}_G \rightarrow A \cotimes \mathcal{K}_G $. 
The maps $ \iota \sigma_n $ are uniformly bounded and converge towards the identity on 
$ A \otimes \mathbb{K}_G $ pointwise. Hence they converge towards the 
identity uniformly on small subsets of $ A \otimes \mathbb{K}_G $.
Explicitly, if $ p_n \in A^+ \cotimes \mathcal{K}_G $ denotes the element
given by 
$$
p_n = \sum_{j = 1}^n 1 \otimes |e_j\ket \bra e_j| 
$$
then $ \sigma_n $ can be written as $ \sigma_n(T) = p_n T p_n $. Since the vectors $ e_j $ are smooth we conclude that
$ \sigma_n $ induces a bounded linear map $ \Smooth(A \otimes \mathbb{K}_G) \rightarrow \Smooth(A) \cotimes \mathcal{K}_G $ which 
will again be denoted by $ \sigma_n $. The maps $ \iota \sigma_n $ converge towards the 
identity uniformly on small subsets of $ \Smooth(A \otimes \mathbb{K}_G) $ as well. 
Hence the claim follows from the isoradial subalgebra theorem \ref{isosubalgebratheorem1}. \qed \\
We conclude this section with another application of the isoradial subalgebra theorem. Recall from \cite{Voigtbs} that a $ G $-simplicial complex is 
a simplicial complex $ X $ with a type-preserving smooth simplicial action of the totally disconnected group $ G $. 
We will assume in the sequel that all $ G $-simplicial complexes have at most countably many simplices. A regular smooth 
function on $ X $ is a function whose restriction to each simplex $ \sigma $ of $ X $ is smooth in the usual 
sense and which is constant in the direction orthogonal to the boundary $ \partial \sigma $ in a neighborhood 
of $ \partial \sigma $. The algebra $ C^\infty_c(X) $ of regular smooth functions on $ X $ with compact support 
is a $ G $-algebra in a natural way. 
\begin{prop}
Let $ X $ be a finite dimensional and locally finite $ G $-simplicial complex. Then the natural map 
$ \iota: C^\infty_c(X) \rightarrow C_0(X) $ is an isoradial $ G $-subalgebra and defines an invertible element 
in 
$$ 
HL^G_0(C^\infty_c(X), \Smooth(C_0(X))).
$$
\end{prop}
\proof As for smooth manifolds one checks that the inclusion homomomorphism $ \iota: C^\infty_c(X) \rightarrow C_0(X) $ 
is isoradial. By induction over the dimension of $ X $ we shall construct a sequence of bounded linear maps 
$ \sigma_n: C_0(X) \rightarrow C^\infty_c(X) $ such that $ \iota \sigma_n $ converges to the identity 
uniformly on small sets. For $ k = 0 $ this is easily achieved by restriction of functions to finite subsets 
and extension by zero. Assume that the maps $ \sigma_n $ are constructed for all $ (k - 1) $-dimensional 
$ G $-simplicial complexes and assume that $ X $ is $ k $-dimensional. If $ X^{k - 1} $ denotes the $ (k - 1) $-skeleton of $ X $
we have a commutative diagram 
\begin{equation*} 
\xymatrix{
C^\infty(X,X^{k - 1}) \;\; \ar@{>->}[r]\ar@{->}[d] & C^\infty_c(X) \ar@{->>}[r]\ar@{->}[d] & C^\infty_c(X^{k - 1}) \ar@{->}[d] \\
C(X,X^{k - 1}) \;\; \ar@{>->}[r] & C_0(X) \ar@{->>}[r] & C_0(X^{k - 1}) 
     }
\end{equation*}
where $ C^\infty(X,X^{k - 1}) $ and 
$ C(X,X^{k - 1}) $ denote the kernels of the canonical restriction homomorphisms and the vertical arrows are 
natural inclusions. It is shown in \cite{Voigtbs} that the upper extension has a bounded linear splitting, and 
the lower extension has a bounded linear splitting as well. Note that the $ C(X,X^{k - 1}) $ is a $ C^* $-direct 
sum of algebras of the form $ C_0(\Delta^k \setminus \partial \Delta^k) $ where $ \Delta^k $ denotes the standard $ k $-simplex 
and $ \partial \Delta^k $ is its boundary. Similarly, $ C^\infty_c(X,X^{k - 1}) $ 
is the bornological direct sum of corresponding subalgebras $ C^\infty_c(\Delta^k \setminus \partial \Delta^k) $. 
Hence, by applying suitable cutoff functions, we are reduced to 
construct approximate inverses to the inclusion $ C^\infty_c(\Delta^k\setminus \partial \Delta^k) \rightarrow C_0(\Delta^k \setminus \partial \Delta^k) $. 
This is easily achieved using smoothing operators. Taking into account the isoradial subalgebra theorem \ref{isosubalgebratheorem1} yields
the assertion. \qed 

\section{Tensor products}\label{secextprod}

In this section we study the equivariant $ X $-complex for the analytic tensor algebra of the tensor product of two $ G $-algebras. 
This will be used in the construction of the equivariant Chern-Connes character in the odd case. \\
Let us first recall the definition of the tensor product of paracomplexes of $ \ayd $-modules \cite{Voigtepch}. 
If $ C $ and $ D $ are paracomplexes of separated $ \ayd $-modules then the tensor product $ C \boxtimes D $ is given by 
$$ 
(C \boxtimes D)_0 = C_0 \otimes_{\mathcal{O}_G} D_0 \oplus C_1 \otimes_{\mathcal{O}_G} D_1, \qquad 
(C \boxtimes D)_1 = C_1 \otimes_{\mathcal{O}_G} D_0 \oplus C_0 \otimes_{\mathcal{O}_G} D_1
$$
where the group $ G $ acts diagonally and $ \mathcal{O}_G $ acts by multiplication. Using that 
$ \mathcal{O}_G $ is commutative one checks that the tensor product $ C \boxtimes D $ becomes a separated $ \ayd $-module in this way.
The boundary operator $ \partial $ in $ C \boxtimes D $ is defined by 
$$
\partial_0 = 
\begin{pmatrix}
\partial \otimes \id & - \id \otimes \partial \\
\id \otimes \partial & \partial \otimes T
\end{pmatrix} 
\qquad 
\partial_1 = 
\begin{pmatrix}
\partial \otimes T & \id \otimes \partial \\
-\id \otimes \partial & \partial \otimes \id
\end{pmatrix}
$$
and turns $ C \boxtimes D $ into a paracomplex. Remark that the formula for $ \partial $ 
does not agree with the usual definition of the differential in a tensor product of complexes. \\
Now let $ A $ and $ B $ be separated bornological algebras. As it is explained in \cite{CQ1}, the unital free 
product $ A^+ * B^+ $ of $ A^+ $ and $ B^+ $ can be written as  
$$
A^+ * B^+ = A^+ \otimes B^+ \oplus \bigoplus_{j > 0} \Omega^j(A) \otimes \Omega^j(B)
$$
with the direct sum bornology and multiplication given by the Fedosov product
$$
(x_1 \otimes y_1) \circ (x_2 \otimes y_2) = x_1 x_2 \otimes y_1 y_2 - (-1)^{|x_1|} x_1 dx_2 \otimes dy_1 y_2. 
$$
An element $ a_0da_1 \cdots da_n \otimes b_0 db_1 \cdots db_n $ corresponds 
to $ a_0 b_0 [a_1, b_1]\cdots [a_n,b_n] $ in the free product under this identification 
where $ [x,y] = xy - yx $ denotes the ordinary commutator. Note that if $ A $ and $ B $ are $ G $-algebras then the free product is again a 
separated $ G $-algebra in a natural way. \\
Consider the extension 
\begin{equation*}
  \xymatrix{
     I\;\; \ar@{>->}[r] & A^+ * B^+ \ar@{->>}[r]^\pi & A^+ \otimes B^+ 
     }
\end{equation*}
where $ I $ is the kernel of the canonical homomorphism $ \pi: A^+ * B^+ \rightarrow A^+ \otimes B^+ $. 
Using the description of the free product in terms of differential forms one has 
$$
I^k = \bigoplus_{j \geq k} \Omega^j(A) \otimes \Omega^j(B) 
$$
for the powers of the ideal $ I $. \\
Analogous to the analytic bornology on tensor algebras we consider an analytic bornology on free products. 
By definition, the analytic bornology on $ A^+ * B^+ $ is the bornology generated by the sets 
\begin{equation*}
S \otimes T \cup \bigcup_{n = 1}^\infty (S(dS)^n \cup (dS)^n) \otimes (T (dT)^n \cup (dT)^n)
\end{equation*}
for all small sets $ S \subset A $ and $ T \subset B $. This bornology turns $ A^+ * B^+ $ into a separated bornological 
algebra. We write $ A^+ \star B^+ $ for the free product of $ A^+ $ and $ B^+ $ equipped with the analytic bornology. 
Clearly the identity map $ A^+ * B^+ \rightarrow A^+ \star B^+ $ 
is a bounded homomorphism. Consequently the natural homomorphisms $ \iota_A: A^+ \rightarrow A^+ \star B^+ $ 
and $ \iota_B: B^+ \rightarrow A^+ \star B^+ $ are bounded. \\
Every unital homomorphism $ f: A^+ \star B^+ \rightarrow C $ into a unital bornological algebra $ C $ 
is determined by a pair of homomorphisms $ f_A: A \rightarrow C $ and $ f_B: B \rightarrow C $. 
Define a linear map $ c_f: A \otimes B \rightarrow C $ by $ c_f(a, b) = [f_A(a), f_B(b)] $. Let us call 
$ f_A $ and $ f_B $ almost commuting if 
$$
c_f(S)^\infty = \bigcup_{n = 1}^\infty c_f(S)^n
$$ 
is small for every small subset $ S \subset A \otimes B $. Clearly, $ c_f = 0 $ iff the images of $ f_A $ and 
$ f_B $ commute. The following property of $ A^+ \star B^+ $ is a direct consequence of the definition of the analytic bornology. 
\begin{lemma} \label{analyticfreeprod}
Let $ A $ and $ B $ be separated bornological algebras. For a pair of bounded equivariant homomorphisms $ f_A: A \rightarrow C $ and 
$ f_B: B \rightarrow C $ into a unital bornological algebra $ C $ 
the corresponding unital homomorphism $ f: A^+ \star B^+ \rightarrow C $ is bounded iff $ f_A $ and $ f_B $ are 
almost commuting. 
\end{lemma} 
In particular, the canonical homomorphism $ \pi: A^+ \star B^+ \rightarrow A^+ \otimes B^+ $ is bounded and we 
obtain a corresponding extension
\begin{equation*}
  \xymatrix{
     I \;\; \ar@{>->}[r] & A^+ \star B^+ \ar@{->>}[r]^\pi & A^+ \otimes B^+ 
     }
\end{equation*}
of bornological algebras with bounded linear splitting. It is straightforward to verify that the ideal $ I $ 
with the induced bornology is analytically nilpotent. Remark that if $ A $ and $ B $ are $ G $-algebras then all the previous constructions are 
compatible with the group action. \\
Let $ I $ be a $ G $-invariant ideal in a separated $ G $-algebra $ R $ and define the paracomplex 
$ \mathcal{H}^2_G(R,I) $ by 
$$ 
\mathcal{H}^2_G(R,I)^0 = \mathcal{O}_G \otimes R/(\mathcal{O}_G \otimes I^2 + b(\mathcal{O}_G \otimes IdR))
$$ 
in degree zero and by  
$$
\mathcal{H}^2_G(R,I)^1 = \mathcal{O}_G \otimes \Omega^1(R)/(b(\Omega^2_G(R)) + \mathcal{O}_G \otimes I \Omega^1(R))
$$
in degree one with boundary operators induced from $ X_G(R) $. \\
Now let $ A $ and $ B $ be separated $ G $-algebras. We abbreviate $ R = A^+ \star B^+ $ and define an $ \ayd $-map $ \phi: X_G(A^+) \boxtimes X_G(B^+) \rightarrow \mathcal{H}^2_G(R,I) $ by 
\begin{align*}
&\phi(f(t) \otimes x \otimes y) = f(t) \otimes xy \\
&\phi(f(t) \otimes x_0dx_1 \otimes y_0dy_1) = f(t) \otimes x_0 (t^{-1} \cdot y_0)[x_1, t^{-1} \cdot y_1]  \\
&\phi(f(t) \otimes x \otimes y_0 dy_1) = f(t) \otimes xy_0 dy_1 \\
&\phi(f(t) \otimes x_0dx_1 \otimes y) = f(t) \otimes x_0dx_1y 
\end{align*}
where $ [x,y] = xy - yx $ denotes the commutator. The following result for the analytic free product 
$ R = A^+ \star B^+ $ is obtained in the same way as the corresponding assertion in \cite{Voigtepch} for the ordinary free product. 
\begin{prop} \label{tensorprop}
The map $ \phi: X_G(A^+) \boxtimes X_G(B^+) \rightarrow \mathcal{H}^2_G(R,I) $ defined above is an  
isomorphism of paracomplexes for all separated $ G $-algebras $ A $ and $ B $. 
\end{prop}
After these preparations we shall prove the following assertion. 
\begin{prop}\label{extprodnu}
Let $ A $ and $ B $ be separated locally multiplicative $ G $-algebras. Then there exists a natural chain map 
$$
X_G(\mathcal{T}(A^+ \otimes B^+))^{\mathbb{L}c} \rightarrow 
(X_G((\mathcal{T}A)^+)^{\mathbb{L}c} \boxtimes X_G((\mathcal{T}B)^+)^{\mathbb{L}c})^{\mathbb{L}c} 
$$
of paracomplexes. There is an analogous chain map if the derived completion is replaced by the ordinary 
completion. 
\end{prop}
\proof Let us abbreviate $ Q = (\mathcal{T}A)^+ \otimes (\mathcal{T}B)^+ $. The canonical homomorphism 
$ \tau: Q \rightarrow A^+ \otimes B^+ $ induces a bounded equivariant homomorphism $ \mathcal{T}Q \rightarrow 
\mathcal{T}(A^+ \otimes B^+) $. Conversely, the obvious splitting for $ \tau $ is a lanilcur since the algebras $ A $ 
and $ B $ are locally multiplicative. 
It follows that there is a canonical bounded equivariant homomorphism $ \mathcal{T}(A^+ \otimes B^+) 
\rightarrow \mathcal{T}Q $ as well. As a consequence we obtain a natural homotopy equivalence 
$$
X_G(\mathcal{T}Q) \simeq X_G(\mathcal{T}(A^+ \otimes B^+))
$$
using homotopy invariance. \\
We have another analytically nilpotent extension of $ Q $ defined as follows. 
Since commutators in the unital 
free product $ (\mathcal{T}A)^+ * (\mathcal{T}B)^+ $ are mapped to zero under the natural 
map $ (\mathcal{T}A)^+ * (\mathcal{T}B)^+ \rightarrow Q $ we have the extension 
\begin{equation*}
  \xymatrix{
     I\;\; \ar@{>->}[r] & R \ar@{->>}[r]^\pi & Q
     }
\end{equation*}
where $ R = (\mathcal{T}A)^+ \star (\mathcal{T}B)^+ $ is the analytic free product of $ (\mathcal{T}A)^+ $ and $ (\mathcal{T}B)^+ $ and 
$ I $ is the kernel of the bounded homomorphism $ \pi: R \rightarrow Q $. Since the $ G $-algebra $ I $ is analytically nilpotent the 
natural equivariant homomorphism $ \mathcal{T}Q \rightarrow R $ is bounded and induces a chain map 
$ X_G(\mathcal{T}Q) \rightarrow X_G(R) $. \\
Next we have an obvious chain map 
$$
p: X_G(R) \rightarrow \mathcal{H}^2_G(R,I)
$$
and by proposition \ref{tensorprop} there exists a natural isomorphism 
$$
X_G((\mathcal{T}A)^+ ) \boxtimes X_G((\mathcal{T}B)^+) \cong \mathcal{H}^2_G(R,I)
$$
of paracomplexes. Assembling these maps and homotopy equivalences yields a chain map 
$ X_G(\mathcal{T}(A^+ \otimes B^+)) \rightarrow X_G((\mathcal{T}A)^+) \boxtimes X_G((\mathcal{T}B)^+) $. 
Inspecting the construction of the derived completion we get in addition a natural chain map
$$
(X_G((\mathcal{T}A)^+) \boxtimes X_G((\mathcal{T}B)^+))^{\mathbb{L}c} \rightarrow 
(X_G((\mathcal{T}A)^+)^{\mathbb{L}c} \boxtimes X_G((\mathcal{T}B)^+)^{\mathbb{L}c})^{\mathbb{L}c}
$$
which immediately yields the assertion for the derived completion. For the ordinary completion the argument is 
essentially the same. \qed 
\begin{cor} \label{exteriorprod} 
Let $ A $ and $ B $ be separated locally multiplicative 
$ G $-algebras. Then there exists a natural chain map  
$$
X_G(\mathcal{T}(A \otimes B))^{\mathbb{L}c} \rightarrow (X_G(\mathcal{T}A)^{\mathbb{L}c} \boxtimes X_G(\mathcal{T}B)^{\mathbb{L}c})^{\mathbb{L}c}.
$$
An analogous assertion holds if the derived completion is replaced by the ordinary completion. 
\end{cor}
\proof The claim follows easily from proposition \ref{extprodnu} by applying the excision 
theorem \ref{Excision2} to tensor products of the extensions $ 0 \rightarrow A \rightarrow A^+ \rightarrow \mathbb{C} 
\rightarrow 0 $ and $ 0 \rightarrow B \rightarrow B^+ \rightarrow \mathbb{C} \rightarrow 0 $. \qed 
\begin{prop} \label{complexext} 
Let $ A $ be a separated locally multiplicative $ G $-algebra. Then the natural chain map 
$$
X_G(\mathcal{T}(\mathbb{C} \otimes A))^{\mathbb{L}c} \rightarrow 
(X_G(\mathcal{T}\mathbb{C})^{\mathbb{L}c} \boxtimes X_G(\mathcal{T}A)^{\mathbb{L}c} )^{\mathbb{L}c}
$$
is a homotopy equivalence. Similarly, one obtains a homotopy equivalence if the derived 
completion is replaced by the ordinary completion.
\end{prop}
\proof Recall that the natural map $ X_G(\mathcal{T}\mathbb{C})^{\mathbb{L}c} \rightarrow 
X_G(\mathcal{T}\mathbb{C})^c $ is a local homotopy equivalence and that $ X_G(\mathcal{T}\mathbb{C})^c \simeq X_G(\mathbb{C}) = \mathcal{O}_G[0] $ 
using the projection homomorphism
$ \mathcal{T}\mathbb{C} \rightarrow \mathbb{C} $. As a consequence we obtain a natural homotopy equivalence 
$ (X_G(\mathcal{T}\mathbb{C})^{\mathbb{L}c} \boxtimes X_G(\mathcal{T}A)^{\mathbb{L}c})^{\mathbb{L}c} \rightarrow 
(\mathcal{O}_G[0] \boxtimes X_G(\mathcal{T}A)^{\mathbb{L}c})^{\mathbb{L}c} $. The composition 
of the latter with the chain map $ X_G(\mathcal{T}(\mathbb{C} \otimes A))^{\mathbb{L}c} \rightarrow 
(X_G(\mathcal{T}\mathbb{C})^{\mathbb{L}c} \boxtimes X_G(\mathcal{T}A)^{\mathbb{L}c})^{\mathbb{L}c} $ obtained 
in corollary \ref{exteriorprod} can be identified with the canonical homotopy equivalence
$ X_G(\mathcal{T}(\mathbb{C} \otimes A))^{\mathbb{L}c} \cong X_G(\mathcal{T}A)^{\mathbb{L}c} \simeq (X_G(\mathcal{T}A)^{\mathbb{L}c})^{\mathbb{L}c} $. 
This proves the claim for the derived completion. For the ordinary completion the argument is analogous. \qed \\
We remark that using the perturbation lemma one may proceed in a similar way 
as for the periodic theory \cite{Voigtepch} in order to construct a candidate for the homotopy inverse to the map
$ X_G(R)^c \rightarrow \mathcal{H}^2_G(R)^c $ induced by the projection $ p $ occuring in the proof of proposition \ref{extprodnu}. 
The problem is that the formula thus obtained does not yield a bounded map in general. 
However, a more refined construction might yield a bounded homotopy inverse. 
For our purposes proposition \ref{complexext} is sufficient. 

\section{Algebraic description of equivariant Kaspararov theory}\label{seckkg}

In this section we review the description of equivariant $ KK $-theory arising from the approach developped by 
Cuntz \cite{Cuntzdocumenta}, \cite{Cuntzbivkcstar}. This approach to $ KK $-theory is based on extensions
and will be used in the definition of the equivariant Chern-Connes character below. \\
One virtue of the framework in \cite{Cuntzdocumenta} is that it allows to construct bivariant versions of $ K $-theory in very 
general circumstances. Moreover, one can adapt the setup to treat equivariant versions of such theories as well. 
The main ingredient in the definition is a class of extensions in the underlying category of algebras
which contains certain fundamental extensions. In particular one needs a suspension 
extension, a Toeplitz extension and a universal extension. In addition one has to specify a tensor product 
which preserves the given class of extensions. \\ 
For equivariant $ KK $-theory the underlying category of algebras is the category $ G\cstaralg $ of separable $ G $-$ C^* $-algebras. 
By definition, morphisms in $ G\cstaralg $ are the equivariant $ * $-homomorphisms. The correct choice 
of extensions is the class $ \mathfrak{E} $ of extensions of $ G $-$ C^* $-algebras with equivariant completely positive 
splitting. As a tensor product one uses the maximal $ C^* $-tensor product. \\
The suspension extension of a $ G $-$ C^* $-algebra $ A $ is 
$$
\xymatrix{
\E_s(A): A(0,1) \;\ar@{>->}[r] & A(0,1] \ar@{->>}[r] & A \\
 } 
$$
where $ A(0,1) $ denotes the tensor product $ A \otimes C_0(0,1) $, and accordingly the algebras 
$ A(0,1] $ and $ A[0,1] $ are defined. The group action on these algebras is given by the pointwise action on $ A $. \\
The Toeplitz extension is defined by 
$$
\xymatrix{
\E_t(A): \mathbb{K} \otimes A \; \ar@{>->}[r] & \mathfrak{T} \otimes A \ar@{->>}[r] & C(S^1) \otimes A \\
 }
$$
where $ \mathfrak{T} $ is the Toeplitz algebra, that is, the universal $ C^* $-algebra 
generated by an isometry. As usual $ \mathbb{K} $ is the algebra of compact operators, and $ \mathbb{K} $ and 
$ \mathfrak{T} $ are equipped with the trivial $ G $-action. \\
Finally, one needs an appropriate universal extension \cite{Cuntzbivkcstar}. Given an algebra $ A $ in $ G\cstaralg $ there exists a tensor 
algebra $ TA $ in $ G \cstaralg $ together with a canonical surjective equivariant $ * $-homomorphism 
$ \tau_A: TA \rightarrow A $ such that the extension 
$$
\xymatrix{
   {\E_u(A): JA \;\;} \ar@{>->}[r] & TA \ar@{->>}[r] & A \\
 } 
$$ 
is contained in $ \mathfrak{E} $ where $ JA $ denotes the kernel of $ \tau_A $. Moreover, this extension is universal in the following sense. 
Given any extension $ \E: 0 \rightarrow K \rightarrow E \rightarrow A \rightarrow 0 $ in $ \mathfrak{E} $ there exists a commutative 
diagram 
$$
\xymatrix{
   {JA \;\;} \ar@{>->}[r] \ar@{->}[d] & TA \ar@{->>}[r] \ar@{->}[d] & A \ar@{=}[d]  \\
{K \;\;} \ar@{>->}[r] & E \ar@{->>}[r] & A 
 } 
$$
The left vertical map $ JA \rightarrow K $ in this diagram is called the 
classifying map of $ \E $. One should not confuse $ TA $ with the analytic tensor algebra used in the 
construction of analytic and local cyclic homology. \\
One defines $ J^2A = J(JA) $ and recursively $ J^nA = J(J^{n - 1}A) $ for $ n \in \mathbb{N} $ as well 
as $ J^0A = A $. Let us denote by $ \phi_A: JA \rightarrow C(S^1) \otimes A $ the 
equivariant $ * $-homomorphism obtained by composing the classifying map $ JA \rightarrow A(0,1) $ of 
the suspension extension with the inclusion map $ A(0,1) \rightarrow C(S^1) \otimes A $ 
given by viewing $ A(0,1) $ as the ideal of functions vanishing in the point $ 1 $. This yields 
an equivariant $ * $-homomorphism $ \epsilon_A: J^2A \rightarrow \mathbb{K} \otimes A $ as the left 
vertical arrow in the commutative diagram 
$$
\xymatrix{
   {J^2A \;\;} \ar@{>->}[r] \ar@{->}[d]^{\epsilon_A} & TJA \ar@{->>}[r] \ar@{->}[d] & JA \ar@{->}[d]^{\phi_A}  \\
{\mathbb{K} \otimes A \;\;} \ar@{>->}[r] & \mathfrak{T} \otimes A \ar@{->>}[r] & C(S^1) \otimes A
 } 
$$
where the bottom row is the Toeplitz extension $ \E_t(A) $. The classifying map $ \epsilon_A $ plays 
an important role in the theory. If $ [A,B]_G $ denotes the set of equivariant homotopy classes of morphisms between $ A $ and $ B $ 
then the previous construction induces a map $ S: [J^k A, \mathbb{K} \otimes B] \rightarrow [J^{k + 2} A, \mathbb{K} \otimes B] $ 
by setting $ S[f] = [(\mathbb{K} \otimes f) \circ \epsilon_{J^{k + 2}A}] $. Here 
one uses the identification $ \mathbb{K} \otimes \mathbb{K} \otimes B \cong \mathbb{K} \otimes B $. \\
We write $ \mathbb{K}_G $ for the algebra of compact operators on the regular representation 
$ L^2(G) $ equipped with its natural $ G $-action. The equivariant stabilization $ A_G $ of 
a $ G $-$ C^* $-algebra $ A $ is defined by $ A_G = A \otimes \mathbb{K} \otimes \mathbb{K}_G $. It 
has the property that $ A_G \otimes \mathbb{K}(\mathcal{H}) \cong A_G $ as $ G $-$ C^* $-algebras for every 
separable $ G $-Hilbert space $ \mathcal{H} $. \\
Using this notation the equivariant bivariant $ K $-group obtained in the approach of Cuntz can be written as 
$$
kk^G_*(A,B) \cong
\varinjlim_j\; [J^{* + 2j}(A_G), \mathbb{K} \otimes B_G] 
$$
where the direct limit is taken using the maps $ S $ defined above. It follows from the results in \cite{Cuntzbivkcstar} 
that $ kk^G_*(A,B) $ is a graded abelian group and that there exists an associative bilinear product 
for $ kk^G_* $. Let us remark that we have inserted 
the algebra $ J^{* + 2j}(\mathbb{K}_G \otimes \mathbb{K} \otimes A) $ in the formula defining 
$ kk^G_* $ instead of $ \mathbb{K}_G \otimes \mathbb{K} \otimes  J^{* + 2j}A $ as in \cite{Cuntzbivkcstar}. Otherwise the 
construction of the product seems to be unclear. \\
We need some more terminology. A functor $ F $ defined on the category 
of $ G $-$ C^* $-algebras with values in an additive category is called 
(continuously) homotopy invariant if $ F(f_0) = F(f_1) $ whenever $ f_0 $ and $ f_1 $ are equivariantly
homotopic $ * $-homomorphisms. It is called $ C^* $-stable if there exists a natural isomorphism 
$ F(A) \cong F(A \otimes \mathbb{K} \otimes \mathbb{K}_G) $ 
for all $ G $-$ C^* $-algebras $ A $. Finally, $ F $ is called split exact if the sequence
$ 0 \rightarrow F(K) \rightarrow F(E) \rightarrow F(Q) \rightarrow 0 $ is split exact
for every extension $ 0 \rightarrow K \rightarrow E \rightarrow Q \rightarrow 0 $ of $ G $-$ C^* $-algebras 
that splits by an equivariant $ * $-homomorpism $ \sigma: Q \rightarrow E $. \\
Equivariant $ KK $-theory \cite{Kasparov2} can be viewed as an additive category $ KK^G $ with separable 
$ G $-$ C^* $-algebras as objects and $ KK^G_0(A,B) $ as the set of morphisms between two objects $ A $ and $ B $. 
Composition of morphisms is given by the Kasparov product. 
There is a canonical functor $ \iota: G\cstaralg \rightarrow KK^G $ which is the identity on 
objects and sends equivariant $ * $-homomorphisms to the corresponding $ KK $-elements. 
Equivariant $ KK $-theory satisfies the following universal property \cite{Thomsen}, \cite{Meyerkkg}. 
\begin{theorem}\label{kkguniversal}
An additive functor $ F $ from $ G $-$ C^*\emph{-\textsf{Alg}} $ into an additive category 
$ \mathcal{C} $ factorizes uniquely over $ KK^G $ iff it is continuously homotopy invariant, $ C^* $-stable and split exact. 
That is, given such a functor $ F $ there exists a unique functor $ \ch_F: KK^G \rightarrow \mathcal{C} $ such 
that $ F = \ch_F \iota $. 
\end{theorem}
It follows from the theory developped in \cite{Cuntzdocumenta} that the functor $ kk^G $ is homotopy invariant, 
$ C^* $-stable and split exact. In fact, it is universal with respect to these properties. As a consequence 
one obtains the following theorem. 
\begin{theorem}
For all separable $ G $-$ C^* $-algebras $ A $ and $ B $ there is a natural isomorphism $ KK^G_*(A,B) \cong kk^G_*(A,B) $.  
\end{theorem}
As already indicated above we will work with the description of equivariant $ KK $-theory provided by $ kk^G_* $ in the sequel. 
In other words, for our purposes we could as well take the definition of $ kk^G_* $ as definition of equivariant $ KK $-theory. 

\section{The equivariant Chern-Connes character}\label{secchernconnes}

In this section we construct the equivariant Chern-Connes character from equivariant $ KK $-theory into 
equivariant local cyclic homology. Moreover we calculate the character in a simple special case. \\
First let us extend the definition of equivariant local cyclic homology $ HL^G_* $ to bornological algebras that are 
equipped with a not necessarily smooth action of the group $ G $. This is done by first applying the smoothing 
functor $ \Smooth $ in order to obtain separated $ G $-algebras. In particular, we may view equivariant local cyclic 
homology as an additive category $ HL^G $ with the same objects as $ G \cstaralg $ and $ HL^G_0(A, B) $ 
as the set of morphisms between two objects $ A $ and $ B $. By construction, there 
is a canonical functor from $ G\cstaralg $ to $ HL^G $. 
\begin{theorem}\label{localhomprop} Let $ G $ be a totally disconnected group. The canonical functor from 
$ G $-$ C^*\emph{-\textsf{Alg}} $ to $ HL^G $ is continuously homotopy invariant, $ C^* $-stable and split exact. 
\end{theorem}
\proof Proposition \ref{homc*} shows together with proposition \ref{homotopyinv} that $ HL^G $ is continuously homotopy invariant. 
We obtain $ C^* $-stability from proposition \ref{stabkg} together with proposition \ref{stability}. Finally, if
$ 0 \rightarrow K \rightarrow E \rightarrow Q \rightarrow 0 $ is a split exact extension 
of $ G $-$ C^* $-algebras then 
$ 0 \rightarrow \Smooth(K) \rightarrow \Smooth(E) \rightarrow \Smooth(Q) \rightarrow 0 $ is a split exact extension of  
$ G $-algebras. Hence split exactness follows from the 
excision theorem \ref{Excision}. \qed \\
Having established this result, the existence of the equivariant Chern-Connes character in the even case is an 
immediate consequence of the universal property of equivariant Kasparov theory. More precisely, according to theorem \ref{localhomprop} 
and theorem \ref{kkguniversal} we obtain an additive map 
\begin{equation*}
\ch^G_0: KK^G_0(A,B) \rightarrow HL^G_0(A, B) 
\end{equation*}
for all separable $ G $-$ C^* $-algebras $ A $ and $ B $. The resulting transformation is multiplicative with respect to 
the Kasparov product and the composition product, respectively. Remark that the equivariant Chern-Connes character 
$ \ch^G_0 $ is determined by the property that it maps $ KK $-elements induced by equivariant 
$ * $-homomorphisms to the corresponding $ HL $-elements. \\
Before we extend this character to a multiplicative transformation on $ KK^G_* $ we shall describe
$ \ch^G_0 $ more concretely using the theory explained in section \ref{seckkg}. Let us fix some notation. 
If $ f: A \rightarrow B $ is an equivariant homomorphism between $ G $-algebras we denote by $ \ch(f) $ the associated class in 
$ H_0(\SHom_G(X_G(\mathcal{T}A)^{\mathbb{L}c}, X_G(\mathcal{T}B)^{\mathbb{L}c})) $. By slight abuse of 
notation we will also write $ \ch(f) $ for the corresponding element in $ HL^G_0(A,B) $. 
Similarly, assume that $ \E: 0 \rightarrow K \rightarrow E \rightarrow Q \rightarrow 0 $ is an extension of $ G $-algebras with 
equivariant bounded linear splitting. We denote by $ \ch(\E) $ the element $ -\delta(\id_K) $ where 
$$ 
\delta: H_0(\SHom_G(X_G(\mathcal{T}K)^{\mathbb{L}c}, X_G(\mathcal{T}K)^{\mathbb{L}c})) \rightarrow 
H_1(\SHom_G(X_G(\mathcal{T}Q)^{\mathbb{L}c}, X_G(\mathcal{T}K)^{\mathbb{L}c})) 
$$ 
is the boundary map in the six-term exact sequence in bivariant homology obtained from the generalized excision 
theorem \ref{Excision2}. Again, by slight abuse of notation we will also write $ \ch(\E) $ for 
the corresponding element in $ HL^G_1(Q,K) $. \\
If $ f: A \rightarrow B $ is an equivariant $ * $-homomorphism between $ G $-$ C^* $-algebras we write simply $ \ch(f) $ 
instead of $ \ch(\Smooth(f)) $ for the element associated to the corresponding homomorphism of $ G $-algebras. 
In a similar way we proceed for extensions of $ G $-$ C^* $-algebras with equivariant 
completely positive splitting. \\
Using theorem \ref{localhomprop} one shows that 
$ \ch(\epsilon_A) \in HL^G_*(J^2A, \mathbb{K} \otimes A) $ is invertible. The same holds true for the iterated 
morphisms $ \ch(\epsilon_A^n) \in HL^G_*(J^{2n} A, \mathbb{K} \otimes A) $. 
Remark also that $ \ch(\iota_A) \in HL^G_*(A, \mathbb{K} \otimes A) $ is invertible. \\
Now assume that $ x \in KK^G_0(A,B) $ is represented by $ f: J^{2n}A_G \rightarrow \mathbb{K} \otimes B_G $. Then the 
class $ \ch^G_0(f) $ is corresponds to 
$$
\ch(\iota_{A_G}) \cdot \ch(\epsilon^n_{A_G})^{-1} \cdot \ch(f) \cdot \ch(\iota_{B_G})^{-1}
$$
in $ HL^G_0(A_G,B_G) $, and the latter group is canonically isomorphic to $ HL^G_*(A,B) $. \\
For the definition of $ \ch^G_1 $ we follow the discussion in \cite{Cuntzdocumenta}. 
We denote by $ j: C_0(0,1) \rightarrow C(S^1) $ the inclusion homomorphism
obtained by viewing elements of $ C_0(0,1) $ as functions on the circle vanishing in $ 1 $. 
Moreover let $ \mathbb{K} $ be the algebra of compact operators on $ l^2(\mathbb{N}) $ and 
let $ \iota: \mathbb{C} \rightarrow \mathbb{K} $ be the homomorphism determined by sending $ 1 $ to the minimal projection 
onto the first basis vector in the canonical orthonormal basis. 
If $ A $ is any $ G $-$ C^* $-algebra we write 
$ \iota_A: A \rightarrow A \otimes \mathbb{K} $ for the homomorphism obtained 
by tensoring $ \iota $ with the identity on $ A $. \\
In the sequel we write $ \E_s $ instead of $ \E_s(\mathbb{C}) $ and similarly $ \E_t $ instead of $ \E_t(\mathbb{C}) $ 
for the Toeplitz extension of $ \mathbb{C} $.  
\begin{prop}\label{2pilemma} With the notation as above one has 
$$ 
\ch(\E_s) \cdot \ch(j) \cdot \ch(\E_t) = \frac{1}{2 \pi i} \ch(\iota) 
$$ 
in $ H_0(\SHom_G(X_G(\mathcal{T}\mathbb{C})^{\mathbb{L}c}, X_G(\mathcal{T}\mathbb{K})^{\mathbb{L}c})) $.  
\end{prop}
\proof First observe that the same argument as in the proof of proposition \ref{stabkg} shows that the 
element $ \ch(\iota) $ is invertible. 
Let us write $ z $ for the element in $ H_0(\SHom_G(X_G(\mathcal{T}\mathbb{C})^{\mathbb{L}c}, X_G(\mathcal{T}\mathbb{C})^{\mathbb{L}c})) $ given by 
of $ (2 \pi i) \ch(\E_s) \cdot \ch(j) \cdot \ch(\E_t) \cdot \ch(\iota)^{-1} $. It suffices to show that $ z $ is equal to the identity. \\
We consider the smooth analogues of the extensions $ \E_s $ and $ \E_t $ used in \cite{Cuntzdocumenta}. The smooth version of the suspension extension 
is 
$$
\xymatrix{
\mathbb{C}^\infty(0,1) \;\ar@{>->}[r] & \mathbb{C}^\infty(0,1] \ar@{->>}[r] & \mathbb{C} \\
 } 
$$
where $ \mathbb{C}^\infty(0,1) $ denotes the algebra of smooth functions on $ [0,1] $ vanishing with all derivatives 
in both endpoints. Similarly, $ \mathbb{C}^\infty(0,1] $ is the algebra of all smooth functions $ f $ vanishing with all derivatives in $ 0 $ 
and vanishing derivatives in $ 1 $, but arbitrary value $ f(1) $. 
The smooth Toeplitz extension is 
$$
\xymatrix{
\mathbb{K}^\infty \;\ar@{>->}[r] & \mathfrak{T}^\infty \ar@{->>}[r] & C^\infty(S^1) \\
 } 
$$
where $ \mathbb{K}^\infty $ is the algebra of smooth compact operators 
and $ \mathfrak{T}^\infty $ is the smooth Toeplitz algebra defined in \cite{Cuntzdocumenta}. 
We obtain another endomorphism $ z^\infty $ of $ X_G(\mathcal{T}\mathbb{C})^{\mathbb{L}c} $ by repeating the 
construction of $ z $ using the smooth supension and Toeplitz extensions. By naturality one has in fact $ z^\infty = z $, hence it suffices to 
show that $ z^\infty $ is equal to the identity. \\
Recall that we have a local homotopy equivalence 
$ X_G(\mathcal{T}\mathbb{C})^{\mathbb{L}c} \rightarrow X_G(\mathcal{T}\mathbb{C})^c \simeq X_G(\mathbb{C}) $. 
Using the fact that the $ G $-action is trivial on all algebras under consideration the same argument 
as in \cite{Cuntzdocumenta} yields that $ z^\infty $ is equal to the identity. \qed \\
We shall use the abbreviation $ x_A = \ch(\E_u(A)) $ for the element arising from the universal extension of 
the $ G $-$ C^* $-algebra $ A $. 
\begin{prop} \label{epsilontheorem} 
Let $ A $ be a $ G $-$ C^* $-algebra and let $ \epsilon_A: J^2(A) \rightarrow \mathbb{K} \otimes A $ be the canonical map. Then we have the 
relation 
$$
x_A \cdot x_{JA} \cdot \ch(\epsilon_A) = \frac{1}{2 \pi i} \ch(\iota_A)
$$
in $ H_0(\SHom_G(X_G(\mathcal{T}(\Smooth(A) \cotimes \mathcal{K}_G))^{\mathbb{L}c}, 
X_G(\mathcal{T}(\Smooth(A \otimes \mathbb{K}) \cotimes \mathcal{K}_G))^{\mathbb{L}c}))$. 
\end{prop}
\proof For an arbitrary $ G $-$ C^* $-algebra $ A $ 
consider the commutative diagram 
$$
\xymatrix{
X_G(\mathcal{T}(\Smooth(A) \cotimes \mathcal{K}_G))^{\mathbb{L}c}\; \ar@{->}[r]^\cong \ar@{->}[d]^{x_A} 
& X_G(\mathcal{T}(\Smooth(\mathbb{C} \otimes A)\cotimes \mathcal{K}_G))^{\mathbb{L}c} \ar@{->}[d] \\
X_G(\mathcal{T}(\Smooth(JA) \cotimes \mathcal{K}_G))^{\mathbb{L}c} \; \ar@{->}[r] \ar@{->}[d]^{x_{JA}}
& X_G(\mathcal{T}(\Smooth(J\mathbb{C} \otimes A)\cotimes \mathcal{K}_G))^{\mathbb{L}c} \ar@{->}[d] \\
X_G(\mathcal{T}(\Smooth(J^2A)\cotimes \mathcal{K}_G))^{\mathbb{L}c}\; \ar@{->}[r] \ar@{->}[d]^{\ch(\epsilon_A)} 
& X_G(\mathcal{T}(\Smooth(J^2\mathbb{C} \otimes A)\cotimes \mathcal{K}_G))^{\mathbb{L}c} \ar@{->}[d]^{\ch(\epsilon_\mathbb{C} \otimes \id)} \\
X_G(\mathcal{T}(\Smooth(\mathbb{K} \otimes A)\cotimes \mathcal{K}_G))^{\mathbb{L}c} \; \ar@{->}[r] \ar@{->}[d]^{\ch(\iota_A)^{-1}} 
& X_G(\mathcal{T}(\Smooth(\mathbb{K} \otimes A)\cotimes \mathcal{K}_G))^{\mathbb{L}c}  \ar@{->}[d]^{\ch(\iota_{\mathbb{C}} \otimes \id)^{-1}} \\
X_G(\mathcal{T}(\Smooth(A)\cotimes \mathcal{K}_G))^{\mathbb{L}c}\; \ar@{->}[r]^\cong 
& X_G(\mathcal{T}(\Smooth(\mathbb{C} \otimes A) \cotimes \mathcal{K}_G))^{\mathbb{L}c} 
 }
$$  
where the upper part is obtain from the morphism of extensions 
$$
\xymatrix{
JA \; \ar@{>->}[r] \ar@{->}[d] & TA \ar@{->>}[r] \ar@{->}[d] & A \ar@{->}[d]^\cong  \\
A \otimes J\mathbb{C} \; \ar@{>->}[r] & A \otimes T\mathbb{C} \ar@{->>}[r] & A \otimes \mathbb{C} 
 } 
$$
and a corresponding diagram with $ A $ replaced by $ JA $. Observe that there is a natural homomorphism
$ D \cotimes \Smooth(A) \cotimes \mathcal{K}_G  \rightarrow \Smooth(D \otimes A) \cotimes \mathcal{K}_G $ 
for every trivial $ G $-$ C^* $-algebra $ D $. For simplicity we will write $ \Smooth(A) $ 
instead of $ \Smooth(A) \cotimes \mathcal{K}_G $ in the following commutative diagram 
$$
\xymatrix{
X_G(\mathcal{T}(\mathbb{C} \cotimes \Smooth(A)))^{\mathbb{L}c}\; \ar@{->}[r] \ar@{->}[d] 
& (X_G(\mathcal{T}\mathbb{C})^{\mathbb{L}c} \boxtimes X_G(\mathcal{T}\Smooth(A))^{\mathbb{L}c})^{\mathbb{L}c} \ar@{->}[d]^{x_{\mathbb{C}}\boxtimes \id} \\
X_G(\mathcal{T}(J\mathbb{C} \cotimes \Smooth(A)))^{\mathbb{L}c} \; \ar@{->}[r] \ar@{->}[d]
& (X_G(\mathcal{T}(J\mathbb{C}))^{\mathbb{L}c} \boxtimes X_G(\mathcal{T}\Smooth(A))^{\mathbb{L}c})^{\mathbb{L}c} \ar@{->}[d]^{x_{J\mathbb{C}} \boxtimes \id} \\
X_G(\mathcal{T}(J^2\mathbb{C}\cotimes \Smooth(A)))^{\mathbb{L}c}\; \ar@{->}[r] \ar@{->}[d]^{\ch(\epsilon_\mathbb{C} \cotimes \id)}
& (X_G(\mathcal{T}(J^2\mathbb{C}))^{\mathbb{L}c} \boxtimes X_G(\mathcal{T}\Smooth(A))^{\mathbb{L}c})^{\mathbb{L}c} 
\ar@{->}[d]^{\ch(\epsilon_\mathbb{C})\boxtimes \id} \\
X_G(\mathcal{T}(\mathbb{K}\cotimes \Smooth(A)))^{\mathbb{L}c} \; \ar@{->}[r] \ar@{->}[d]^{\ch(\iota_A)^{-1}} 
& (X_G(\mathcal{T}\mathbb{K})^{\mathbb{L}c} \boxtimes X_G(\mathcal{T}\Smooth(A))^{\mathbb{L}c})^{\mathbb{L}c} \ar@{->}[d]^{\ch(\iota)^{-1} \boxtimes \id} \\
X_G(\mathcal{T}(\mathbb{C}\cotimes \Smooth(A)))^{\mathbb{L}c}\; \ar@{->}[r] 
& (X_G(\mathcal{T}\mathbb{C})^{\mathbb{L}c} \boxtimes X_G(\mathcal{T}\Smooth(A))^{\mathbb{L}c})^{\mathbb{L}c}
 }
$$  
obtained using corollary \ref{exteriorprod}. According to proposition \ref{complexext} the first and the last horizontal map 
in this diagram are homotopy equivalences. Moreover, we may connect the right column of the first diagram with the 
left column of the previous diagram. Using these observations the assertion follows from proposition \ref{2pilemma} in the same way as 
in \cite{Cuntzdocumenta}. \qed \\
After these preparations we shall now define the Chern-Connes character in the odd case. 
For notational simplicity we assume that all $ G $-$ C^* $-algebras $ A $ are replaced by their equivariant 
stabilizations $ A_G $. We may then use the identification 
$$
KK^G_*(A,B) \cong \varinjlim_j\; [J^{* + 2j}(A), \mathbb{K} \otimes B] 
$$
and obtain a canonical isomorphism $ KK^G_1(A,B) \cong KK^G_0(JA,B) $. Consider an element $ u \in KK^G_1(A,B) $ and denote by 
$ u_0 $ the element in $ KK^G_0(JA,B) $ corresponding to $ u $. Then the element $ \ch^G_1(u) \in HL^G_1(A,B) $ is defined by  
$$
\ch^G_1(u) = \sqrt{2 \pi i}\; x_A \cdot \ch^G_0(u_0)
$$
in terms of the character in the even case obtained before. Using proposition \ref{epsilontheorem} one concludes in the same way 
as in \cite{Cuntzdocumenta} that the formula
$$
\ch^G_{i + j}(x \cdot y) = \ch^G_i(x) \cdot \ch^G_j(y)
$$
holds for all elements $ x \in KK^G_i(A,B) $ and $ y \in KK^G_j(B,C) $. \\
We have now completed the construction of the equivariant Chern-Connes character and summarize the result in 
the following theorem. 
\begin{theorem}
Let $ G $ be a second countable totally disconnected locally compact group and let $ A $ and $ B $ be 
separable $ G $-$ C^* $-algebras. Then there exists a transformation  
\begin{equation*}
\ch^G_*: KK^G_*(A,B) \rightarrow HL^G_*(A,B) 
\end{equation*}
which is multiplicative with respect to the Kasparov product in $ KK^G_* $ and the composition product 
in $ HL^G_* $. Under this transformation elements in $ KK^G_0(A,B) $ induced by equivariant $ * $-homomorphisms from $ A $ to $ B $ are mapped to the 
corresponding elements in $ HL^G_0(A,B) $. 
\end{theorem}
The transformation obtained in this way will be called the equivariant Chern-Connes character. 
One shows as in nonequivariant case that, up to possibly a sign and a factor $ \sqrt{2 \pi i} $, 
the equivariant Chern-Connes character is compatible with the boundary maps 
in the six-term exact sequences associated to an extension in $ \mathfrak{E} $. \\
At this point it is not clear wether the equivariant Chern-Connes character is a useful tool to detect 
information contained in equivariant $ KK $-theory. As a matter of fact, equivariant local cyclic homology groups 
are not easy to calculate in general. In a separate paper we will exhibit interesting situations in which $ \ch^G_* $ becomes in 
fact an isomorphism after tensoring the left hand 
side with the complex numbers. At the same time a convenient description of the right 
hand side of the character will be obtained. \\
Here we shall at least illustrate the nontriviality of the equivariant Chern-Connes character in a simple special case. 
Assume that $ G $ be a profinite group. The character of a finite dimensional representation of $ G $ defines an element in 
the algebra $ \mathcal{R}(G) = (\mathcal{O}_G)^G $ of conjugation invariant smooth functions on $ G $. As usual we 
denote by $ R(G) $ the representation ring of $ G $.
\begin{prop} Let $ G $ be a profinite group. Then the equivariant Chern-Connes character 
\begin{equation*}
\ch^G_*: KK^G_*(\mathbb{C},\mathbb{C}) \rightarrow HL^G_*(\mathbb{C}, \mathbb{C}) 
\end{equation*}
can be identified with the character map $ R(G) \rightarrow \mathcal{R}(G) $. This identification is compatible with the products. 
\end{prop}
\proof Let $ V $ be a finite dimensional representation of $ G $. Then $ \mathbb{K}(V) $ is a unital $ G $-algebra 
and the element in $ R(G) = KK^G_0(\mathbb{C},\mathbb{C}) $ 
corresponding to $ V $ is given by the class of the equivariant homomorphism $ p_V: \mathbb{C} \rightarrow \mathbb{K}(V) $ 
in $ KK^G_0(\mathbb{C}, \mathbb{K}(V)) \cong KK^G_0(\mathbb{C},\mathbb{C}) $ where $ p_V $ is defined by 
$ p_V(1) = \id_V $. Using stability of $ HL^G_* $ proposition \ref{hlhacompcomp} we see that the class of 
$ \ch^G_0(p_V) $ in $ HL^G_0(\mathbb{C}, \mathbb{K}(V)) \cong 
HL^G_0(\mathbb{C}, \mathbb{C}) = \mathcal{R}(G) $ corresponds to the character of 
the representation $ V $. The claim follows easily from these observations. \qed 

\bibliographystyle{plain}

\end{document}